\documentclass[11pt,a4paper,twoside]{amsart}

\usepackage[headings]{fullpage}

\usepackage{amsmath}
\usepackage{amsthm}
\usepackage{amsfonts}
\usepackage{amssymb}
\usepackage{mathtools}
\usepackage[english]{babel}
\usepackage{combelow}
\usepackage{multirow}
\usepackage{stackengine}
\usepackage[hyphens]{url}

\newtheorem{theorem}{Theorem}[section]
\newtheorem{lemma}[theorem]{Lemma}

\theoremstyle{definition}

\newcommand{\Q}{{\mathbb{Q}}}
\newcommand{\R}{{\mathbb{R}}}

\newcommand{\floor}[1]{\left\lfloor#1\right\rfloor}

\newcommand{\prpr}{{\prime\prime}}
\newcommand{\bp}{{b^\prime}}
\newcommand{\bpp}{{b^\prpr}}
\newcommand{\bl}{{\tilde{b}}}

\newcommand{\si}{{\sqrt{-1}}}
\newcommand{\I}{{\mathcal{I}}}
\DeclareMathOperator*{\h}{h}
\DeclareMathOperator*{\Imf}{Im}
\newcommand{\Mod}[1]{\ \left(\mathrm{mod}\ #1\right)}

\let\lm=\lambda
\let\Lm=\Lambda

\let\abs=\envert
\let\sq=\sqrt

\newcommand{\acr}{\newline\indent}

\theoremstyle{remark}

\begin{document}
\title{Three-term Machin-type formulae}
\author[Tomohiro Yamada]{Tomohiro Yamada*}
\address{\llap{*\,}Center for Japanese language and culture\acr
                   Osaka University\acr
                   562-8678\acr
                   3-5-10, Semba-Higashi, Minoo, Osaka\acr
                   JAPAN}
\email{tyamada1093@gmail.com}


\begin{abstract}
We shall show that there exist only finitely many nondegenerate three-term Machin-type
formulae and give explicit upper bounds for the sizes of variables.
\end{abstract}

\maketitle

\section{Introduction}\label{intro}

Machin's formula
\begin{equation}\label{eq01}
4\arctan\frac{1}{5}-\arctan\frac{1}{239}=\frac{\pi}{4}
\end{equation}
is well known and has been used to calculate approximate values of $\pi$.
Analogous formulae
\begin{equation}\label{eq02}
\arctan\frac{1}{2}+\arctan\frac{1}{3}=\frac{\pi}{4},
\end{equation}
\begin{equation}\label{eq03}
2\arctan\frac{1}{2}-\arctan\frac{1}{7}=\frac{\pi}{4}
\end{equation}
and
\begin{equation}\label{eq04}
2\arctan\frac{1}{3}+\arctan\frac{1}{7}=\frac{\pi}{4},
\end{equation}
which are also well known,
were attributed to Euler, Hutton and Hermann, respectively.
But according to Tweddle \cite{Twe}, these formulae also
seem to have been found by Machin.

Several three-term formulae such as
\begin{equation}
8\arctan\frac{1}{10}-\arctan\frac{1}{239}-4\arctan\frac{1}{515}=\frac{\pi}{4}
\end{equation}
due to Simson in 1723 (see \cite{Twe}) and
\begin{equation}
12\arctan\frac{1}{18}+8\arctan\frac{1}{57}-5\arctan\frac{1}{239}=\frac{\pi}{4}
\end{equation}
due to Gauss in 1863 also have been known.

More generally, an $n$-term Machin-type formula is defined to be an identity of the form
\begin{equation}\label{eq10}
y_1 \arctan \frac{1}{x_1}+y_2 \arctan \frac{1}{x_2}+\cdots +y_n \arctan \frac{1}{x_n}=\frac{r\pi}{4}
\end{equation}
with integers $x_1, x_2, \ldots, x_n, y_1, y_2, \ldots, y_n$ and $r\neq 0$.

Theoretical studies of Machin-type formulae have begun with a series of works of St{\o}rmer,
who proved that \eqref{eq01}-\eqref{eq04} are all two-term ones in 1895 \cite{Stm1}
and gave a necessary and sufficient condition for given integers $x_1, x_2, \ldots, x_n>1$
to have a Machin-type formula \eqref{eq10} and $102$ three-term ones in 1896 \cite{Stm2}.
St{\o}rmer asked for other three-term Machin-type formulae and questioned whether
there exist infinitely many ones or not.
Up to now the only known other nontrivial (i.e. not derived from \eqref{eq02}-\eqref{eq04}) three-term formulae are
\begin{equation}
5\arctan\frac{1}{2}+2\arctan\frac{1}{53}+\arctan\frac{1}{4443}=\frac{3\pi}{4},
\end{equation}
\begin{equation}
5\arctan\frac{1}{3}-2\arctan\frac{1}{53}-\arctan\frac{1}{4443}=\frac{\pi}{2},
\end{equation}
and
\begin{equation}
5\arctan\frac{1}{7}+4\arctan\frac{1}{53}+2\arctan\frac{1}{4443}=\frac{\pi}{4}.
\end{equation}
\cite{WH} attributes these formulae to Wrench \cite{Wre} although these formulae cannot be found there.
We note that the second and the third formulae follow from the first formula using \eqref{eq02} and \eqref{eq03}
respectively.

The purpose of this paper is to answer to St{\o}rmer's other question in negative.
We shall show that there exist only finitely many three-term Machin-type formulae
which do not arise from a linear combinations of three identities \eqref{eq02}-\eqref{eq04}.

St{\o}rmer's criterion is essentially as follows:
For given integers $x_1, x_2, \ldots, x_n>1$, \eqref{eq10} holds for some integers $y_1, y_2, \ldots, y_n$ and $r\neq 0$
if and only if there exist integers $s_{i, j} (i=1, 2, \ldots, n, j=1, 2, \ldots, n-1)$
and Gaussian integers $\eta_1, \eta_2, \ldots, \eta_{n-1}$ such that
\begin{equation}\label{eq120}
\left[\frac{x_i+\si}{x_i-\si}\right]=\left[\frac{\eta_1}{\bar\eta_1}\right]^{\pm s_{i, 1}}\left[\frac{\eta_2}{\bar\eta_2}\right]^{\pm s_{i, 2}} \cdots \left[\frac{\eta_{n-1}}{\bar\eta_{n-1}}\right]^{\pm s_{i, n-1}}.
\end{equation}
for $i=1, 2, \ldots, n$.

Writing $m_j=\eta_j \bar\eta_j$ for $j=1, 2, \ldots, n-1$,
this condition can be reformulated as follows:
there exist nonnegative integers $s_{i, j}  (i=1, 2, \ldots, n, j=1, 2, \ldots, n)$
with $0\leq s_{i, n}\leq 1$ such that the equation
\begin{equation}\label{eq130}
x_i^2+1=2^{s_{i, n}} m_1^{s_{i, 1}} m_2^{s_{i, 2}}\cdots m_{n-1}^{s_{i, n-1}}
\end{equation}
holds for $i=1, 2, \ldots n$ and, additionally,
$x_i\equiv \pm x_j\Mod{m_k}$ for three indices $i, j, k$
with $x_i^2+1\equiv x_j^2+1\equiv 0\Mod{m_k}$.

Thus, for given three integers $x_1, x_2, x_3>1$,
there exist nonzero integers $y_1$, $y_2$, $y_3$, and $r$
such that $\gcd(y_1, y_2, y_3)=1$ and a three-term Machin-type formula
\begin{equation}\label{eq11}
y_1 \arctan \frac{1}{x_1}+y_2 \arctan \frac{1}{x_2}+y_3 \arctan \frac{1}{x_3}=\frac{r\pi}{4}
\end{equation}
holds if and only if
there exist integers $k_i, \ell_i (i=1, 2, 3)$ and Gaussian integers $\eta_1, \eta_2$ such that
\begin{equation}\label{eq12}
\left[\frac{x_i+\si}{x_i-\si}\right]=\left[\frac{\eta_1}{\bar\eta_1}\right]^{\pm k_i}\left[\frac{\eta_2}{\bar\eta_2}\right]^{\pm \ell_i}
\end{equation}
holds for $i=1, 2, 3$ or, equivalently,
writing $m_j=\eta_j \bar\eta_j$ for $j=1, 2, \ldots, n-1$ and choosing $v_i\in \{0, 1\}$ appropriately,
the equation
\begin{equation}\label{eq13}
x_i^2+1=2^{v_i} m_1^{k_i} m_2^{\ell_i}
\end{equation}
holds for $i=1, 2, 3$ and, additionally,
$x_i\equiv \pm x_{i^\prime} \Mod{m_j}$ for two indices $i, i^\prime$ with
$x_i^2+1\equiv x_{i^\prime}^2+1\equiv 0\Mod{m_j}$.
Furthermore, there exists a nonzero integers $d$ such that
$dy_1=\pm k_2 \ell_3\pm k_3 \ell_2$, $dy_2=\pm k_3 \ell_1\pm k_1 \ell_3$,
and $dy_3=\pm k_1 \ell_2\pm k_2 \ell_1$
holds simultaneously with appropriate choices of signs.

Now we shall state our result in more detail.
\begin{theorem}\label{thm1}
Assume that $x_1, x_2, x_3, y_1, y_2, y_3$ and $r$ are nonzero integers
with $x_1, x_2, x_3>1$, $\{x_1, x_2, x_3\}\neq \{2, 3, 7\}$ satisfying \eqref{eq11}, and
$\gcd(y_1, y_2, y_3)=1$
and $m_1, m_2, s_i, k_i, \ell_i (i=1, 2, 3)$ are corresponding integers with $m_2>m_1>0$
satisfying \eqref{eq13}.
\begin{itemize}
\item[I.] If $m_1 m_2\geq 10^{10}$ and $m_2$ divides $x_i^2+1$ for $i=1, 2, 3$, then
$m_1<m_2<1.342\times 10^{34}, x_i<\exp(8380114833)$, and $\abs{y_i}<9.152\times 10^{16}$.
\item[II.] If $m_1 m_2\geq 10^{10}$ and $m_2$ does not divide $x_i^2+1$ for some $i$, then
$m_1<2.531\times 10^{24}$, $m_2<\exp(294622)$,
$x_i<\exp(5.054\times 10^{12})$, and $\abs{y_i}<1.312\times 10^{19}$.
Moreover, $m_2<e^{40000}$ if additionally $m_1\geq 9.134\times 10^{22}$.
\item[III.] If $m_1 m_2<10^{10}$, then
$x_i<\exp(8813999998)$ and $\abs{y_i}<4.508\times 10^{18}$.
\end{itemize}
More detailed estimates are given in Tables \ref{tbl5} and \ref{tbl6}.
\end{theorem}

We use a lower bound for linear forms in three logarithms in order to obtain
upper bounds for exponents $k_i$'s and $\ell_i$'s in terms of $m_1, m_2$.

These upper bounds themselves do not give finiteness of $m_1$ and $m_2$.
However, noting that $r\neq 0$, which gives $\abs{\sum_i y_i \arctan(1/x_i)}\geq \pi/4$,
the first case can be easily settled using these upper bounds.
In order to settle the second case, we additionally need an upper bound $m_2$ in terms of $m_1$.
This can be done using a lower bound for a quantity of the form $y\arctan(1/x)-r\pi/2$,
which gives a linear form of two logarithms.
Computation was done by PAGI-GP.
Out PARI-GP script is available from \\
\url{https://drive.google.com/drive/folders/1gFTjX05ts5eaz8ncOjfaxRMvCW-drvGX}.

\eqref{eq130} can be seen as a special case of the generalized Ramanujan-Nagell equation
\begin{equation}\label{eq14}
x^2+Ax+B=p_1^{e_1} p_2^{e_2} \cdots p_n^{e_n},
\end{equation}
where $A$ and $B$ are given integers with $A^2-4B\neq 0$ and $p_1, p_2, \ldots, p_n$ are given primes.
Evertse \cite{Eve} proved that \eqref{eq14} has at most $3\times 7^{4n+6}$ solutions.
In the case $n=2$, the author \cite{Ymd} reduced Evertse's bound $3\times 7^{14}$ to $63$.

On the other hand, our result does not give an upper bound for numbers of solutions
\begin{equation}\label{eq15}
x^2+1=2^s p_1^k p_2^\ell
\end{equation}
since the case $r=0$ is not considered.
Indeed, St{\o}rmer \cite{Stm2} implicitly pointed out that, if $x^2+1=ay$, then
\begin{equation}
\arctan \frac{1}{az-x}-\arctan \frac{1}{az+a-x}=\arctan \frac{1}{az(z+1)-(2z+1)x+y}.
\end{equation}

St{\o}rmer \cite{Stm3} showed that
\eqref{eq15} has at most one solution
with each fixed combination of parities of $s_i, k_i, \ell_i$ with zero and nonzero-even distinguished.
Although there exist $18$ combinations $(0\mid 1, 0\mid 1\mid 2, 0\mid 1\mid 2)$,
all-even combinations can clearly be excluded and therefore \eqref{eq15} has at most
$14$ solutions totally.

\section{Preliminaries}

In this section, we introduce some notation and some basic facts.

For any Gaussian integer $\eta$, we have an associate $\eta^\prime$ of $\eta$ such that
$-\pi/4<\arg \eta^\prime<\pi/4$ and therefore $-\pi/2<\arg \eta^\prime/\bar{\eta^\prime}<\pi/2$.
As usual, the absolute logarithmic height $\h(\alpha)$ of an algebraic number $\alpha$
of degree $d$ over $\Q$ is defined by
\begin{equation}
\h(\alpha)=\frac{1}{d}\left(\log\abs{a}+\sum_{i=1}^d\max\{0, \log\abs{\alpha^{(i)}}\}\right),
\end{equation}
where $a$ is the smallest integer such that $a\alpha$ is an algebraic integer and
the $\alpha^{(i)}$'s $(i=1, 2, \ldots, d)$ are the complex conjugates of $\alpha$ over $\Q$.

We call a formula \eqref{eq10} to be degenerate if
\begin{equation}
\sum_{i\in S}y_i^\prime \arctan\frac{1}{x_i}=\frac{r^\prime \pi}{4}
\end{equation}
for some proper subset $S$ of $\{1, 2, \ldots, n\}$
and integers $y_i^\prime ~ (i\in S)$ and $r^\prime$ which may be zero but not all zero.

We shall show that a degenerate case of a three-term Machin-type formula occurs
only when $\{x_1, x_2, x_3\}=\{2, 3, 7\}$.
\begin{lemma}\label{lm21}
Let $x_1, x_2, x_3, y_1, y_2, y_3$ and $r$ be nonzero integers with $x_1, x_2, x_3>1$ satisfying \eqref{eq11}
and $m_1, m_2, s_i, k_i, \ell_i (i=1, 2, 3)$ be corresponding integers satisfying \eqref{eq13}.
If $\ell_1=\ell_2=0$, then $\{x_1, x_2, x_3\}=\{2, 3, 7\}$.
\end{lemma}

\begin{proof}
Let $\eta_1, \eta_2$ be corresponding Gaussian integers in \eqref{eq12}.
Since $\ell_1=\ell_2=0$, we have $[(x_i+\si)/(x_i-\si)]=[\eta_1/\bar \eta_1]^{\pm k_i}$
for $i=1, 2$.
Hence, we obtain
\begin{equation}
k_2\arctan\frac{1}{x_1}\pm k_1\arctan\frac{1}{x_2}=\frac{r^\prime \pi}{4}
\end{equation}
for some integer $r^\prime$, which may be zero.
Substituting this into \eqref{eq11} multiplied by an appropriate factor, we have
\begin{equation}
y_1^\prime \arctan\frac{1}{x_1}+y_3^\prime \arctan\frac{1}{x_3}=\frac{r^\prpr \pi}{4}
\end{equation}
for some integers $y_1^\prime, y_3^\prime, r^\prpr$, which may be zero again but not all zero.
Now St{\o}rmer's result in \cite{Stm1} on two-term Machin-type formulae implies that
$\{x_1, x_2, x_3\}=\{2, 3, 7\}$.
\end{proof}

Moreover, we can determine the only solution for \eqref{eq12} with $k_i=0$ and $\ell_i>2$.
\begin{lemma}\label{lm22}
The equation $x^2+1=2^e y^n$ with $x>0, n>2$ has
only one integral solution $(x, e, y, n)=(239, 1, 13, 4)$.
\end{lemma}
\begin{proof}
From the result of \cite{Leb}, $x^2+1=m^t$ with $x>0, t>1$ has no solution.
Th\'{e}or\`{e}me 8 of \cite{Stm3} shows that $x^2+1=2m^t$, then $t$ must be a power of two.
By Ljunggren's result \cite{Lju}, the only integer solution of $x^2+1=2m^4$
with $x, m>1$ is $(x, m)=(239, 13)$.
Easier proofs of Ljunggren's result have been obtained by Steiner and Tzanakis \cite{ST}
and Wolfskill \cite{Wol}.
\end{proof}

\section{A lower bound for linear forms of three logarithms}

Our argument depends on a lower bound for linear forms of three logarithms.
Results in a paper \cite{MV2} entitled \textit{a kit} are technical but still worthful to use for the purpose of improving our upper bounds.
We use Theorem 4.1 of an earlier version \cite{MV1} rather than \cite{MV2}, which enables us to obtain better bounds.
However, the proof of Theorem 4.1 in \cite{MV1} needs an amendment.

\begin{lemma}\label{lm31}
Let $\alpha_1$, $\alpha_2$, and $\alpha_3$ be three distinct algebraic numbers not equal to
zero or one, which are either all $>1$ or all complex of modulus $1$.
Let $b_1$, $b_2$, and $b_3$ be three positive rational integers with $\gcd(b_1, b_2, b_3)=1$,
and put
$$\Lm=b_1\log\alpha_1+b_2\log\alpha_2-b_3\log\alpha_3,$$
where the logarithms of the $\alpha_i$'s are real if $\alpha_i>1$ and purely imaginary if $\abs{\alpha_i}=1$.
Let $w$ is the maximal order of a root of unity which belongs to the field
$\Q(\alpha_1, \alpha_2, \alpha_3)$ and assume that $0<\abs{\Lm}<2\pi/w$.

We put $d_i=\gcd(b_i, b_3)$ and $b_i=d_i b_i^\prime$ for $i=1, 2$ and $b_3=d_1 b_3^\prime=d_2 b_3^\prpr$.
Let $K\geq 3$, $L\geq 5$, and $R_i$ and $R_{ij}$ with $1\leq i, j\leq 3$ be
positive rational integers satisfying
$$R_i>R_{i1}+R_{i2}+R_{i3}$$
for $i=1, 2, 3$.
Let $\rho\geq 2$ be a real number and suppose that
\begin{equation}\label{eq31}
\left(KL-K-\frac{K}{3L}\right)\log\rho\geq (D+1)\log N+gL(a_1 R_1+a_2 R_2+a_3 R_3)+D(K-1)\log b,
\end{equation}
where $N=K^2 L$, $D=[\Q(\alpha_1, \alpha_2, \alpha_3)\colon\Q]/[\R(\alpha_1, \alpha_2, \alpha_3)\colon\R]$,
$$g=\frac{1}{4}-\frac{N}{12R_1 R_2 R_3},
~b=(b_3^\prime \eta_0)(b_3^\prpr \zeta_0)\left(\prod_{k=1}^{K-1}k!\right)^{-4/(K^2-K)}$$
with
$$\eta_0=\frac{R_1-1}{2}+\frac{(R_3-1)b_1}{2b_3},~\zeta_0=\frac{R_2-1}{2}+\frac{(R_3-1)b_2}{2b_3},$$
and
$$a_i\geq \rho\abs{\log\alpha_i}-\log\abs{\alpha_i}+2D\h(\alpha_i)$$
for $i=1, 2, 3$.
Put
$$V=\sq{(R_{11}+1)(R_{21}+1)(R_{31}+1)}.$$

If, for some positive real number $\chi$,
\begin{itemize}
\item[(i)] $(R_{11}+1)(R_{21}+1)(R_{31}+1)>KM$,
where
$$M=\max\{R_{11}+R_{21}+1, R_{21}+R_{31}+1, R_{31}+R_{11}+1, \chi V\},$$
\item[(ii)] $\# \{\prod_{i=1}^3 \alpha_i^{r_i}: 0\leq r_i\leq R_{i1}\}>L$,
\item[(iii)] $(R_{12}+1)(R_{22}+1)(R_{32}+1)>2K^2$,
\item[(iv)] $\# \{\prod_{i=1}^3 \alpha_i^{r_i}: 0\leq r_i\leq R_{i2}\}>2KL$, and
\item[(v)] $(R_{13}+1)(R_{23}+1)(R_{33}+1)>6K^2 L$,
\end{itemize}
then either
\begin{itemize}
\item[(A)] $\Lm^\prime>\rho^{-KL}$, where
$$\Lm^\prime=\abs{\Lm}\frac{LR_3 e^{LR_3 \abs{\Lm}/(2b_3)}}{2b_3},$$
\item[(B)] $b_i\leq\max\{R_{i1}, R_{i2}\}$ for each $i=1, 2, 3$, or
\item[(C)] there exist rational integers $u_i$, not all zero, such that
$$u_1 b_1+u_2 b_2+u_3 b_3=0$$
with $\gcd(u_1, u_2, u_3)=1$ and
$$\abs{u_i}\leq\frac{(R_{i+1, 1}+1)(R_{i+2, 1}+1)}{M-\max\{R_{i+1, 1}, R_{i+2, 1}\}}$$
for each $i=1, 2, 3$, where the indices are reduced modulo $3$.
\end{itemize}
\end{lemma}

\begin{proof}
We use notation such as $p_i(z_1, z_2)$, $\Psi_\I(x)$, and $\Theta(D, \I)$ according to \cite{MV1}.
In p.8, l.-7 of \cite{MV1}, they wrote that
for a certain linear transformation $T$, the identity
$$T(\theta_1 z_1+\theta_2 z_2)=\theta_1(z_1 +z_2)+\theta_2(-\theta_1 z_1 / \theta_2)=\theta_1 z_1$$
holds, where the right hand side $\theta_1 z_1$ should be $\theta_1 z_2$.
Thus, $p_i(z_1, z_2)$ cannot be reduced to a polynomial $p_i^{(T)}(z_1, z_2)$ with degree at most $D_1$ in $z_2$ in their way.

However, the degree of $p_i^{(T)}(z_1, z_2)$ is at most $\max\{D_1, D_2\}$
and, Lemma 3.5 in \cite{MV1} holds with $\Theta(\max\{D_1, D_2\}, \# \I)$ in place of $\max\{\Theta(D_1, \# \I), \Theta(D_2, \# \I)\}$.
This suffices for our purpose since we set $D_1=D_2=K-1$.

More explicitly, if $\theta_1\theta_2\neq 0$, then we set
$\eta_j^\prime=\eta_j+\theta_2\zeta_j/\theta_1$ and $\zeta_j^\prime=-\theta_2\eta_j/\theta_1$,
and
$q_i(z_1, z_2)=p_i(z_1+z_2, -\theta_1 z_2/\theta_2)$.
Then we have $q_i$ are polynomials with degree $\max\{D_1, D_2\}$ in $z_2$ such that $q_i(x\eta_j^\prime, x\zeta_j^\prime)=p_i(x\eta_j, x\zeta_j)$
for $i=1, \ldots, N$.
Hence, for $1\leq i, j\leq N$, we have
\begin{equation}\label{eq301}
p_i(x\eta_j, x\zeta_j)f_i(x(\theta_1\eta_j+\theta_2\zeta_j)=q_i(x\eta_j^\prime, x\zeta_j^\prime)f_i(\theta x\eta_j^\prime)
\end{equation}
with $\theta=\theta_1$ and each $q_i$ has degree at most $\max\{D_1, D_2\}$ in $z_2$.

It may be better to give an explicit argument in the case $\theta_1 \theta_2=0$.
If $\theta_1=0$, then, for each $1\leq i, j\leq N$, putting $\eta_j^\prime=\zeta_j$, $\zeta_j^\prime=\eta_j$, and $q_i(z_1, z_2)=p_i(z_2, z_1)$,
we have \eqref{eq301} with $\theta=\theta_2$ and each $q_i$ has degree at most $D_1$ in $z_2$.
Similarly, if $\theta_2=0$, then, for each $1\leq i, j\leq N$, putting $\eta_j^\prime=\eta_j$, $\zeta_j^\prime=\zeta_j$, and $q_i(z_1, z_2)=p_i(z_1, z_2)$,
we have \eqref{eq301} with $\theta=\theta_1$ and each $q_i$ has degree at most $D_2$ in $z_2$.

Now, \eqref{eq301} yields that
\begin{equation}\label{eq302}
\Psi_\I(x)=\det(\delta_{i, j}q_i(x\eta_j^\prime, x\zeta_j^\prime)f_i(\theta x\eta_j^\prime)
\end{equation}
and $\Psi_\I(x)$ has a zero at $x=0$ of multiplicity at least $\Theta(\max\{D_1, D_2\}, \# \I)$.

For the remaining part of the proof, we can proceed as in \cite{MV2}
with changes of constants according to \cite{MV1} and using Laurent's zero lemma, which is given as Proposition A.1 in \cite{MV2} in the place of Proposition 3.11 in \cite{MV2}.
\end{proof}

We also need a lower bound for linear forms in two logarithms.
\begin{lemma}\label{lm32}
Let $h, \rho$ and $\mu$ be real numbers with $\rho>1, 1/3\leq \mu\leq 1$.
Set
\[
\sigma=\frac{1+2\mu-\mu^2}{2}, \lm=\sigma\log\rho, H=\frac{h}{\lm}+\frac{1}{\sigma},
\omega=2\left(1+\sq{1+\frac{1}{4H^2}}\right), \theta=\frac{1}{2H}+\sq{1+\frac{1}{4H^2}}.
\]
Let
\begin{equation}
\Lm_1=b_1\log \alpha_1+b_2\log \alpha_2,
\end{equation}
where $b_1, b_2$ are positive integers and
$\alpha_1=(X_1+Y_1\si)/(X_1-Y_1\si), \alpha_2=(X_2+Y_2\si)/(X_2-Y_2\si)$ are Gaussian rationals
with $X_i>\abs{Y_i}>0$ integers.
Put $N_i=X_i^2+Y_i^2, a_i^\prime=\rho_1\abs{\log \alpha_i}+2\h(\alpha_i)=\rho_1 \abs{\log \alpha_i}+\log N_i$
(we use this notation here in order to avoid confusion with $a_i$) and
\begin{equation}
\bpp=\frac{\abs{b_1}}{a_2^\prime}+\frac{\abs{b_2}}{a_1^\prime}.
\end{equation}
Furthermore, assume that $a_1^\prime a_2^\prime\geq \lm^2$ and $h\geq \log \bpp+\log \lm+1.81$.

Then we have
\begin{equation}\label{eq32}
\log\abs{\Lm_1}>-c\left(h+\frac{\lm}{\sigma}\right)^2 a_1^\prime a_2^\prime-\sq{\omega\theta}\left(h+\frac{\lm}{\sigma}\right)
-\log\left(c^\prime \left(h+\frac{\lm}{\sigma}\right)^2a_1^\prime a_2^\prime\right),
\end{equation}
where
\begin{equation}
c=\frac{\mu}{\lm^3\sigma}
\left(\frac{\omega}{6}+\sq{\frac{\omega^2}{36}
+\frac{2\lm\omega^{5/4}\theta^{1/4}}{3\sq{a_1^\prime a_2^\prime H}}
+\frac{\lm\omega}{3H}\left(\frac{1}{a_1^\prime}+\frac{1}{a_2^\prime}\right)}\right)^2,
c^\prime=\sq{\frac{c\sigma\omega\theta}{\lm^3\mu}}.
\end{equation}
\end{lemma}
\begin{proof}
We apply Theorem 2 of \cite{Lau} with $D=1$ and $h=\log \bpp+\log \lm+1.81$.
We note that our settings always yield that
$\log \bpp+\log \lm+1.81>\max\{\lm, (\log 2)/2\}$ in Theorem 2 of \cite{Lau}.
Thus, we obtain \eqref{eq32}.
\end{proof}

Now we are to describe our explicit bounds.

\begin{theorem}\label{thm31}
Let $B_0=e^{20}$ and $X_1, Y_1, X_2, Y_2$ be nonzero rational integers such that
$m_i=X_i^2+Y_i^2$ are odd, $\gcd(X_i, Y_i)=1$ for both $i=1, 2$,
and neither of $m_1$ nor $m_2$ is a perfect power.
Moreover assume that $m_1\leq 41$ and $m_2\geq \tilde m_2$ with $\tilde m_2$ from Table \ref{tbl1}
or $m_i\geq\tilde m_i$ for $i=1, 2$ with $\tilde m_1$ and $\tilde m_2$ from Table \ref{tbl1}.
Write $\eta_i=(X_i+Y_i\si)/(X_i-Y_i\si)$ for $i=1, 2$ and $\eta_3=\si$ and put
\begin{equation}
\Lm=\pm e_1\log\eta_1\pm e_2\log\eta_2\pm e_3\log\eta_3.
\end{equation}

Writing $\Omega=(\log m_1)(\log m_2)$ and $E=(e_1\log m_1+e_2\log m_2)/\Omega$,
we have either $E<A_0$, $(A_2 E)^2\log m_1<B_0$,
\begin{equation}\label{eq33}
-\log\abs{\Lm}<A_1\log^2 (A_2 E\sqrt{\log m_1})\Omega,
\end{equation}
or
\begin{equation}\label{eq34}
e_1 \log m_1+e_2 \log m_2 < A_3 \Omega\log(A_2 E\sqrt{\log m_1}),
\end{equation}
where $A_1, A_2, A_3$, and $\rho$ are constants given in Table \ref{tbl1}.
Here and hereafter, we write $m_1^*=\max\{m_1, 10^{10}/m_1\}$.
\end{theorem}

\begin{table}
\caption{Constants in bounds for Theorem \ref{thm31}}
\begin{center}
\begin{small}
\begin{tabular}{| c | c | c | c | c | c | c |}
\hline
$m_1$ & $\tilde m_2$ & $A_1$ & $A_2$ & $A_3$ & $A_0$ & $\rho$ \\
\hline
$5$ & $2000000005$ & $6792.639$ & $1.987228$ & $356.1620$ & $113.2663$ & $4.531$ \\
$13$ & $769230773$ & $5831.233$ & $1.933988$ & $265.7853$ & $84.8351$ & $4.635$ \\
$17$ & $588235297$ & $3300.665$ & $1.735844$ & $275.7772$ & $75.8539$ & $5.155$ \\
$29$ & $344827589$ & $3812.426$ & $1.820629$ & $242.0317$ & $71.5043$ & $4.931$ \\
$37$ & $270270277$ & $2561.555$ & $1.583766$ & $258.3560$ & $63.9328$ & $5.661$ \\
$41$ & $243902441$ & $5072.693$ & $1.890607$ & $226.4842$ & $69.5800$ & $4.770$ \\
$5$ & $10^{16}$ & $6067.262$ & $1.757445$ & $359.7087$ & $107.0678$ & $4.791$ \\
$13$ & $10^{16}$ & $5143.961$ & $1.702538$ & $274.8311$ & $78.7827$ & $4.932$ \\
$17$ & $10^{16}$ & $2865.651$ & $1.500984$ & $290.9531$ & $69.8349$ & $5.572$ \\
$29$ & $10^{16}$ & $3295.501$ & $1.566799$ & $252.8907$ & $65.0122$ & $5.336$ \\
$37$ & $10^{16}$ & $2176.447$ & $1.410974$ & $258.1061$ & $59.7321$ & $5.878$ \\
$41$ & $10^{16}$ & $4394.894$ & $1.620517$ & $235.0623$ & $62.7626$ & $5.169$ \\
$5$ & $10^{24}$ & $5719.809$ & $1.703384$ & $356.2457$ & $105.9955$ & $4.799$ \\
$13$ & $10^{24}$ & $4843.038$ & $1.649921$ & $271.1685$ & $77.7108$ & $4.938$ \\
$17$ & $10^{24}$ & $2679.811$ & $1.418838$ & $283.1656$ & $68.1770$ & $5.660$ \\
$29$ & $10^{24}$ & $3089.378$ & $1.482171$ & $245.5286$ & $63.3267$ & $5.421$ \\
$37$ & $10^{24}$ & $2028.535$ & $1.331976$ & $250.9353$ & $58.1814$ & $5.975$ \\
$41$ & $10^{24}$ & $4129.103$ & $1.569853$ & $231.1519$ & $61.7086$ & $5.172$ \\
$5$ & $e^{40000}$ & $4994.390$ & $1.479858$ & $375.2430$ & $100.3593$ & $5.135$ \\
$13$ & $e^{40000}$ & $4214.314$ & $1.431915$ & $271.9920$ & $72.8793$ & $5.249$ \\
$17$ & $e^{40000}$ & $2287.357$ & $1.189845$ & $307.5275$ & $62.8590$ & $6.243$ \\
$29$ & $e^{40000}$ & $2652.249$ & $1.242648$ & $264.3578$ & $57.7106$ & $5.986$ \\
$37$ & $e^{40000}$ & $1701.715$ & $1.069216$ & $297.9336$ & $51.8596$ & $6.940$ \\
$41$ & $e^{40000}$ & $3591.493$ & $1.289024$ & $252.7466$ & $55.1890$ & $5.834$ \\\hline
$\tilde m_1$ & $\tilde m_2$ & $A_1$ & $A_2$ & $A_3$ & $A_0$ & $\rho$ \\
\hline
$53$ & \addstackgap[.5\dimexpr 2pt \relax]{$m_1^*$} &
$5381.263$ & $1.985977$ & $219.4190$ & $69.6043$ & $4.589$ \\
$53$ & $10^{24}$ & $4370.434$ & $1.638498$ & $214.9439$ & $61.5175$ & $4.965$ \\
$53$ & $10^{32}$ & $4230.482$ & $1.546752$ & $221.6378$ & $59.1862$ & $5.166$ \\
$53$ & $10^{34}$ & $4205.618$ & $1.541857$ & $221.2940$ & $59.0913$ & $5.167$ \\
$53$ & $e^{2000}$ & $3811.172$ & $1.370913$ & $221.6077$ & $55.4326$ & $5.473$ \\
$53$ & $e^{40000}$ & $3794.956$ & $1.398933$ & $224.5476$ & $55.9594$ & $5.397$ \\
$10^{24}$ & $10^{24}$ & $1469.356$ & $1.105031$ & $139.9162$ & $24.0442$ & $7.292$ \\
$5000$ & $e^{2000}$ & $2375.090$ & $1.189591$ & $177.7400$ & $36.3032$ & $6.256$ \\
$5000$ & $e^{40000}$ & $2363.495$ & $1.186121$ & $177.4782$ & $36.2441$ & $6.259$ \\
$10^7$ & $e^{40000}$ & $1742.083$ & $1.057810$ & $154.1620$ & $27.6471$ & $6.957$ \\
$10^{14}$ & $e^{40000}$ & $1368.425$ & $0.943348$ & $146.8828$ & $21.8950$ & $7.836$ \\

\hline
\end{tabular}
\label{tbl1}
\end{small}
\end{center}
\end{table}

\begin{proof}
We may assume that $X_i>Y_i>0$ for $i=1, 2$ and $\abs{\Lm}<10^{-9}$.
We can take $\alpha_1=\pm\si$, $\alpha_2=\eta_1^{\pm 1}$, and $\alpha_3=\eta_2^{\pm 1}$
with a suitable choice of the combination of signs such that
\begin{equation}
\pm \Lm=e_3\log\alpha_1+e_1\log\alpha_2-e_2\log\alpha_3
\end{equation}
with $(b_1, b_2, b_3)=(e_3, e_1, e_2)$.

Now we set $\rho$ according to Table \ref{tbl1}
and take $a_1=\rho\pi/2$ and $a_3=\rho\pi/2 + \log m_2$.
Moreover, we can take $a_2=2\rho \arctan (Y_1/X_1)+\log m_1$ if $m_1\leq 41$ and
$a_2=\rho\pi/2 + \log m_1$ otherwise.
Hence, we have $\rho>4$, $a_1\geq 2\pi$, $a_2\geq 6\arctan(1/4)+\log 17=4.79304\cdots$.
Moreover, we see that $a_2 a_3\geq (9.062\arctan(1/4)+\log 37)(2.2655\pi+\log 270270277)=135.5131\cdots$.
We note that in the case $(\tilde m_1, \tilde m_2)=(53, m_1^*)$, we have
$a_2 a_3\geq (\log(53)+\rho\pi/2)(\log(10^{10}/53)+\rho\pi/2)$
since $(\log x + \rho\pi/2)(\log(10^{10})-\log x+\rho\pi/2)$ is increasing
in the range $1\leq x\leq 10^5$.

We set $\chi$, $C_0$, and $C_1$ according to Table \ref{tbl2}, $C_2=6^{1/3}$
and put
\begin{equation}
c_1^\prime=(\chi_1 C_1)^{2/3}+\frac{C_2 C_0^{1/3} C_1^{2/3}}{\log^{1/3}\rho}
+\frac{2^{1/3}C_1^{2/3}}{\log^{1/3} B_0}
\end{equation}
and
\begin{equation}
c_i^\prime=(\chi_1 C_1)^{2/3}+\frac{C_2 C_0^{1/3} C_1^{2/3}}{\log^{1/3}\rho}+
\max\left\{\frac{2^{1/3}C_1^{2/3}}{\log^{1/3} B_0}, \sq{\frac{C_0 C_1}{a_1\log\rho}}\right\}
\end{equation}
for $i=2, 3$. 
Moreover, we set
\begin{equation}\label{eqdfbp}
\bp=\left(\frac{c_1 b_3^\prime}{a_1}+\frac{c_3 b_1^\prime}{a_3}\right)
\left(\frac{c_2 b_3^\prpr}{a_2}+\frac{c_3 b_2^\prime}{a_3}\right)
\end{equation}
and $b_0=(A_2 E)^2\log m_1$.
Finally, we put
$K=\floor{C_1 a_1 a_2 a_3 \log b_0}$, $L=\floor{C_0 \log b_0/\log \rho}$,
and $\chi=\chi_1 \log^{1/2} b_0$.
We note that $K\geq K_0$ and $L\geq L_0$ for constants $K_0$ and $L_0$ in Table \ref{tbl2}
in each case.

We put
$R_{i1}=\floor{(\chi^2 K^2 a_{i+1} a_{i+2}/a_i^2)^{1/3}}$
for $i=1, 2, 3$.
Since $\chi_1>0.6$, $C_1>6$, $K>540000$, $C_0\leq 3.5$, and $a_1 a_2 a_3>20$ by our settings
and $b_0\geq B_0=e^{20}$ by assumption,
we have $((1-1/K)\chi_1 C_1)^{2/3}>C_0^{1/2}/(2\sq{a_1 a_2 a_3 \log b_0})$ and therefore
\begin{equation}
R_{i, 1}+1>(\chi^2 K^2 a_{i+1} a_{i+2}/a_i^2)^{1/3}>\sq{\frac{C_0 a_1 a_{5-i}\log b_0}{4a_i}}
\end{equation}
for $i=1, 2$.
Hence, we obtain $4(R_{2, 1}+1)(R_{3, 1}+1)>L$ to confirm (ii).

Observing that $\chi_1>0.6$, $C_1>6$, $K>540000$, $C_0\leq 3.5$,
$a_i>2\pi$ for $i=1, 2, 3$, and $b_0\geq B_0=e^{20}$, we have
\begin{equation}
R_{i1}+R_{i+1, 1}+1<1+(\chi_1 C_1)^{2/3} a_{i+2} (a_i+a_{i+1}) \log b_0
<\chi_1 C_1 a_1 a_2 a_3\log b_0
\end{equation}
for $i=1, 2, 3$ and therefore
\begin{equation}
\begin{split}
(R_{11}+1)(R_{21}+1)(R_{31}+1)> & ~ \chi_1^2 K^2
>K \chi_1^2 C_1 a_1 a_2 a_3 \log^2 b_0 \\
> & ~ K\max\{R_{11}+R_{21}+1, R_{21}+R_{31}+1, R_{31}+R_{11}+1\}.
\end{split}
\end{equation}
Observing that $(R_{11}+1)(R_{21}+1)(R_{31}+1)>(\chi K)^2$, we obtain $(R_{11}+1)(R_{21}+1)(R_{31}+1)>\chi KV$.
Hence, we can confirm (i).

Nextly, we put
\[\begin{split}
R_{12}= & ~ \floor{(2K^2 a_2 a_3/a_1^2)^{1/3}}, \\
R_{22}= & ~ \floor{\max\left\{\sq{\frac{KL a_3}{2a_2}},
\left(\frac{2K^2 a_1 a_3}{a_2^2}\right)^{1/3}\right\}},
R_{32}=\floor{\max\left\{\sq{\frac{KL a_2}{2a_3}},
\left(\frac{2K^2 a_1 a_2}{a_3^2}\right)^{1/3}\right\}}.
\end{split}\]
We immediately see that $4(R_{22}+1)(R_{32}+1)>2KL$ and $(R_{12}+1)(R_{22}+1)(R_{32}+1)>2K^2$
to confirm (iii) and (iv) respectively.

Finally, we put $R_{i3}=\floor{C_2(K^2 L a_{i+1}a_{i+2}/a_i^2)^{1/3}}$ for $i=1, 2, 3$.
We immediately have (v) since $C_2\geq 6^{1/3}$.

In each case, noting that $a_{k+2}>a_{k+1}>a_k$, we have
$R_{k, j}\geq R_{k+1, j}\geq R_{k+2, j}$ for $j=1, 3$
and $R_{k+2, j}$ can be minorized by replacing $m_2$ by $\tilde m_2$ and then $m_1$ by $\tilde m_1$.
In the case $(\tilde m_1, \tilde m_2)=(53, m_1^*)$, $R_{k+2, j}$ can be minorized by
replacing $m_2$ by $10^{10}/m_1$ and then taking $m_1=53$.
Moreover, we observe that $R_{i2}\geq \floor{(2K^2 a_k a_i/a_j^2)^{1/3}}$,
which can be minorized by the same way.
Hence, we have $R_1+R_2+R_3\geq R_0$ for $i=1, 2, 3$ for a constant $R_0$ in Table \ref{tbl2}
in each case.

Now we have
\begin{equation}
\sum_{i=1}^3 a_i R_{i1}<3(\chi_1 C_1)^{2/3} a_1 a_2 a_3 \log b_0,
\end{equation}
either
\begin{equation}
\sum_{i=1}^3 a_i R_{i2}<a_1 a_2 a_3\left(2 \sq{\frac{2C_0 C_1}{\rho\pi\log\rho}} \log b_0
+2^{1/3} C_1^{2/3} \log^{2/3} b_0\right)
\end{equation}
or
\begin{equation}
\sum_{i=1}^3 a_i R_{i2}<3\times 2^{1/3} C_1^{2/3} a_1 a_2 a_3 \log^{2/3} b_0,
\end{equation}
and
\begin{equation}
\sum_{i=1}^3 a_i R_{i3}<3 C_2 C_0^{1/3} C_1^{2/3} a_1 a_2 a_3 \log b_0/\log^{1/3}\rho.
\end{equation}
Taking $R_i=R_{i1}+R_{i2}+R_{i3}+1$ and writing $c_4=(1+1/R_0)^3$, we see that
\begin{equation}
\begin{split}
& R_1 R_2 R_3=\prod_{i=1}^3(R_{i1}+R_{i2}+R_{i3}+1) \\
< & ~ c_4((\chi K)^{2/3}+C_2 K^{2/3} L^{1/3}+\sq{KL} a_1^{-1/3} a_2^{1/6} a_3^{1/6})^2 \\
& ((\chi K)^{2/3}+C_2 K^{2/3} L^{1/3}+(2K^2)^{1/3}) \\
< & ~ c_4\left(C_2+\frac{\chi_1^{2/3}\log^{1/3} b_0}{L^{1/3}}
+\left(\frac{L}{K}\right)^{1/6}a_1^{-1/3} a_2^{1/6} a_3^{1/6}\right)^2 \\
& \left(C_2+\frac{2^{1/3}+\chi_1^{2/3}\log^{1/3} b_0}{L^{1/3}}\right)K^2 L \\
< & ~ c_4\left(C_2+\left(\frac{L\log\rho}{(L-1)C_0}\right)^{1/3}\chi_1^{2/3}
+\left(\frac{K C_0}{C_1(K-1)\log\rho}\right)^{1/6}a_1^{-1/2}\right)^2 \\
& \left(C_2+\left(\frac{L\log\rho}{(L-1)C_0}\right)^{1/3}\chi_1^{2/3}+\left(\frac{2}{L}\right)^{1/3}\right)K^2 L
\end{split}
\end{equation}
or
\begin{equation}
\begin{split}
R_1 R_2 R_3<c_4^3
\left(C_2+\left(\frac{L\log\rho}{(L-1)C_0}\right)^{1/3}\chi_1^{2/3}+\left(\frac{2}{L}\right)^{1/3}\right)^3K^2 L.
\end{split}
\end{equation}
Writing $R_1 R_2 R_3<c_0 K^2 L$, we have $g<1/4-1/(12c_0)\leq g_0$ with $g_0$ taken from Table \ref{tbl3}.

Thus, we conclude that
\begin{equation}
gL \sum_{i=1}^3 a_i R_i<g_0L(a_1+a_2+a_3+C_3^\prime a_1 a_2 a_3 \log b_0)<C_3 a_1 a_2 a_3 \log^2 b_0,
\end{equation}
where
\begin{equation}
C_3^\prime = 3(\chi_1 C_1)^{2/3}+\frac{3C_2 C_0^{1/3} C_1^{2/3}}{\log^{1/3} \rho}+
\max \left\{2 \sq{\frac{2C_0 C_1}{\rho\pi\log\rho}}+\frac{2^{1/3} C_1^{2/3}}{\log^{1/3} B_0},
\frac{3\times 2^{1/3}C_1^{2/3}}{\log^{1/3} B_0}\right\}
\end{equation}
and
\begin{equation}
C_3=\frac{g_0 C_0}{\log\rho} \left(\frac{a_1+a_2+a_3}{a_1 a_2 a_3 \log^2 B_0} + C_3^\prime\right)
\end{equation}
with the note that $C_3\leq C_{3, 0}$ in Table \ref{tbl3}. 

Now we would like to show that $\log b\leq \log b_0$.
We see that $R_i-1=R_{i1}+R_{i2}+R_{i3}<c_i^\prime a_{i+1} a_{i+2} \log b_0$ for $i=1, 2, 3$.
Hence, we obtain
\begin{equation}
b_3^\prime \eta_0=\frac{(R_1 -1)b_3^\prime}{2}+\frac{(R_3 -1)b_1 b_3^\prime}{2b_3}
<\frac{\log b_0}{2}(c_1^\prime a_2 a_3 b_3^\prime+c_3^\prime a_1 a_2 b_1^\prime).
\end{equation}
and
\begin{equation}
b_3^\prpr \zeta_0=\frac{(R_2 -1)b_3^\prpr}{2}+\frac{(R_3 -1)b_2 b_3^\prpr}{2b_3}
<\frac{\log b_0}{2}(c_2^\prime a_3 a_1 b_3^\prpr+c_3^\prime a_1 a_2 b_2^\prpr).
\end{equation}
Thus, we have
\begin{equation}
\begin{split}
\log b\leq & ~ \log\left(\frac{\log^2 b_0}{4}\right)
+\log (c_1^\prime a_2 a_3 b_3^\prime+c_3^\prime a_1 a_2 b_1^\prime)
(c_2^\prime a_3 a_1 b_3^\prpr+c_3^\prime a_1 a_2 b_2^\prpr) \\
& ~ -\frac{4}{K(K-1)}\log\left(\prod_{k=1}^{K-1}k!\right).
\end{split}
\end{equation}
We see that the last term is greater than
\begin{equation}
2\log K-3+\frac{2\log(2\pi K/e^{3/2})}{K-1}-\frac{2+6\pi^{-2}+\log K}{3K(K-1)}>2\log (K+1)-3,
\end{equation}
proceeding as in the proof of Lemme 8 of \cite{LMN} and observing that
\begin{equation}
\frac{2+6\pi^{-2}+\log K}{3K(K-1)}+\frac{1-2\log(2\pi K/e^{3/2})}{K-1}<0
\end{equation}
for $K\geq 2$. 
Since $\log b_0\geq \log B_0=20$ and $C_1>e^{3/2}/2$ by assumption, we have
\begin{equation}\label{eq35}
\begin{split}
\log b\leq & \log\left(\frac{\log^2 b_0}{4}\right)
+\log\left((c_1^\prime a_2 a_3 b_3^\prime+c_3^\prime a_1 a_2 b_1^\prime)
(c_2^\prime a_3 a_1 b_3^\prpr+c_3^\prime a_1 a_2 b_2^\prpr)\right) \\
& -2\log (a_1 a_2 a_3\log b_0)+3-2\log C_1 \\
< & \log\left[\left(\frac{c_1^\prime b_3^\prime}{a_1}+\frac{c_3^\prime b_1^\prime}{a_3}\right)
\left(\frac{c_2^\prime b_3^\prpr}{a_2}+\frac{c_3^\prime b_2^\prime}{a_3}\right)\right]+3-2\log(2C_1)
<\log\bp.
\end{split}
\end{equation}

We see that $c_i^\prime<A_2^{\prime 2}$ for $i=1, 2, 3$ and therefore $\bp<A_2^{\prime 2} \bl$, where
$A_2^\prime$ is taken from Table \ref{tbl3} and
\begin{equation}
\bl=\left(\frac{b_3}{a_1}+\frac{b_1}{a_3}\right)\left(\frac{b_3}{a_2}+\frac{b_2}{a_3}\right).
\end{equation}
Observing that $(b_2, b_3)=(e_1, e_2)$ and $b_1=e_3\leq (2/\pi)(e_1\theta_1+e_2\theta_2)+\abs{\Lm}<e_1+e_2+10^{-9}$, we have
\begin{equation}
\begin{split}
\bl=
& \left(\frac{2e_2}{\rho\pi}+\frac{e_3}{\log m_2}\right)\left(\frac{e_2}{\log m_1}+\frac{e_1}{\log m_2}\right) \\
< & \left(\frac{2e_2}{\rho\pi}+\frac{e_1+e_2+10^{-9}}{\log m_2}\right)
\left(\frac{e_2}{\log m_1}+\frac{e_1}{\log m_2}\right) \\
< & (1+10^{-9})\left(\frac{e_2}{\log m_1}+\frac{e_1}{\log m_2}\right)
\left(e_2\left(\frac{2}{\rho\pi}+\frac{1}{\log m_2}\right)+\frac{e_1}{\log m_2}\right) \\
< & \psi^2(m_1, m_2)\left(\frac{e_2}{\log m_1}+\frac{e_1}{\log m_2}\right)^2,
\end{split}
\end{equation}
where
\begin{equation}
\psi=\psi(m_1, m_2)=
\sq{(1+10^{-9})\max\left\{\frac{2\log m_1}{\rho\pi}+\frac{\log m_1}{\log m_2}, 1\right\}}.
\end{equation}

Hence, we have $\bl<\psi^2(m_1, m_2) E^2$ and, combining with \eqref{eq35},
we see that $\bp<A_2^{\prime 2} \bl<(A_2^\prime \psi E)^2$.
Observing that
\begin{equation}
(A_2^\prime \psi)^2<
A_2^{\prime 2}(1+10^{-9})\left(\frac{2}{\rho\pi}+\frac{1}{\log \tilde m_2}\right) \log m_1
<A_2^2 \log m_1,
\end{equation}
we can confirm that $\bp<(A_2 E)^2\log m_1=b_0$.

Now, observing that $D=1$, it suffices to confirm that
\begin{equation}\label{eq31b}
\left(1-\frac{1}{K}\right)\left(1-\frac{1}{L}-\frac{1}{3L^2}\right)C_0 C_1
\geq \frac{2\log (C_0 (C_1 a_1 a_2 a_3)^2 \log^3 b_0/\log\rho)}{a_1 a_2 a_3 \log^2 b_0}+C_1+C_{3, 0}
\end{equation}
instead of \eqref{eq31}. 
Indeed, we see that the left hand side of \eqref{eq31b} is greater than $Y_0$ while
the right hand side of \eqref{eq31b} is smaller than $X_0$, where we take $X_0$ and $Y_0$
from Table \ref{tbl3}.
Since $X_0<Y_0$ in any case, we confirm \eqref{eq31b}.
So that, either of (A), (B), or (C) in Lemma \ref{lm31} holds.

In the case (A), we have
\begin{equation}
\log\abs{\Lm^\prime}>-KL\log\rho>-C_0 C_1 a_1 a_2 a_3 \log^2 b_0.
\end{equation}
Now we shall show that $-\log\abs{\Lm}<(1+10^{-50})KL\log\rho$.
We may assume that $-\log\abs{\Lm}\geq KL\log\rho$.
Since $KL>10^7$ and $\rho>3$, we have
\begin{equation}
\frac{LR_3 \abs{\Lm}}{2b_3}<5KL\times 3^{-KL}<10^{-1000} KL
\end{equation}
and
\begin{equation}
\log\frac{LR_3}{2b_3}<\log \frac{LR_3}{2}<10^{-1000} KL.
\end{equation}
Hence, we obtain
\begin{equation}
-\log\abs{\Lm} <(1+10^{-50})KL\log\rho<(1+10^{-50})C_0 C_1a_1 a_2 a_3\log^2 b_0
<\frac{A_1}{4}\Omega\log^2 b_0.
\end{equation}
Recalling that $b_0=(A_2 E)^2\log m_1$,
we have $-\log\abs{\Lm}<A_1 \log^2 (A_2 E\sqrt{m_1})\Omega$ to obtain
\eqref{eq33}, where we note that, in the case $(\tilde m_1, \tilde m_2)=(53, m_1^*)$,
\begin{equation}
\frac{a_1 a_2}{\Omega}=\left(1+\frac{\rho\pi}{2\log m_1}\right)
\left(1+\frac{\rho\pi}{2\log m_2}\right)
<\left(1+\frac{3.975\pi}{2\log 53}\right)\left(1+\frac{3.975\pi}{2\log (10^{10}/53)}\right).
\end{equation}

In the case (B), we have
\begin{equation}
b_i \leq \max\left\{(\chi_1 C_1)^{2/3} a_{i+1} a_{i+2} \log b_0,
\sq{\frac{C_0 C_1 a_k}{\log\rho}} a_j \log b_0, (2C_1^2 \log^2 b_0)^{1/3} a_{i+1} a_{i+2}\right\}
\end{equation}
for $i\neq j\neq k\neq i$.
For $e_1=b_i$, observing that $a_{i+1}\leq\log m_2+\rho\pi/2$ and $a_{i+2}=\rho\pi/2$, we have
\begin{equation}
e_1\leq
\max\left\{\frac{(\chi_1 C_1)^{2/3}\rho\pi}{2}, \sq{\frac{C_0 C_1 \rho\pi}{4\log\rho}},
\frac{2^{1/3}C_1^{2/3}}{\log^{1/3} B_0}\right\}
\left(\log m_2 +\frac{\rho\pi}{2}\right)\log b_0.
\end{equation}
Similarly, we have
\begin{equation}
e_2\leq
\max\left\{\frac{(\chi_1 C_1)^{2/3}\rho\pi}{2}, \sq{\frac{C_0 C_1 \rho\pi}{4\log\rho}},
\frac{2^{1/3}C_1^{2/3}}{\log^{1/3} B_0}\right\}
\left(\log m_1 +\frac{\rho\pi}{2}\right)\log b_0.
\end{equation}
It follows that
\begin{equation}\label{eq34a}
\begin{split}
e_1\log m_1+e_2\log m_2< & ~
\max\left\{\frac{(\chi_1 C_1)^{2/3}\rho\pi}{2}, \sq{\frac{C_0 C_1 \pi}{4\log\rho}},
\frac{2^{1/3}C_1^{2/3}}{\log^{1/3} B_0}\right\} \\
& ~ \times \left(2+\frac{\rho\pi}{2\log m_1}+\frac{\rho\pi}{2\log m_2}\right)\Omega\log b_0
\end{split}
\end{equation}
to obtain \eqref{eq34}.

Finally, we must settle the case (C).
We see that
\begin{equation}
(R_{i+1, 1}+1)(R_{i+2, 1}+1)<(1+\epsilon_{i1})(\chi_1 C_1)^{4/3} a_i^2 a_{i+1} a_{i+2} \log^2 b_0
\end{equation}
for $i=1, 2, 3$,
where $\epsilon_{i1}=1/R_{i+1, 1}+1/R_{i+2, 1}+1/(R_{i+1, 1}R_{i+2, 1})$.
Since
\begin{equation}
M\geq \chi V>\chi^2 K>\chi_1^2\left(1-\frac{1}{K}\right) (C_1 a_1 a_2 a_3 \log^2 b_0)
\end{equation}
and
\begin{equation}
\max\{R_{i+1, 1}, R_{i+2, 1}\}\leq (\chi_1 C_1)^{2/3} a_i \max\{a_{i+1}, a_{i+2}\} \log b_0,
\end{equation}
we have
\begin{equation}
M-\max\{R_{i+1, 1}, R_{i+2, 1}\}> \chi_1^2 \left(1-\epsilon_{i2}-\frac{1}{K}\right) C_1 a_1 a_2 a_3 \log^2 b_0,
\end{equation} 
where
\begin{equation}
\epsilon_{i2}=\frac{\chi_1^{2/3}}{C_1^{1/3}\min\{a_{i+1}, a_{i+2}\}\log B_0}.
\end{equation}
Hence, we obtain
\begin{equation}
\abs{u_i}\leq\frac{(R_{i+1, 1}+1)(R_{i+2, 1}+1)}{M-\max\{R_{i+1, 1}, R_{i+2, 1}\}}
<(1+\epsilon)\chi_1^{-2/3} C_1^{1/3} a_i
\end{equation}
for each $i=1, 2, 3$, where $\epsilon=\max_{1\leq i\leq 3}(1+\epsilon_{i1})/(1-\epsilon_{i2}-1/K)-1$.

We begin with the case $u_1\neq 0$.
We use
\begin{equation}
\begin{split}
u_1 \Lm= & ~ u_1(b_1\log\alpha_1+b_2\log\alpha_2-b_3\log\alpha_3) \\
= & ~ b_2(u_1\log\alpha_2-u_2\log\alpha_1)-b_3(u_3\log\alpha_1+u_1\log\alpha_3) \\
= & ~ b_2\log\alpha_1^\prime + b_3\log\alpha_2^\prime,
\end{split}
\end{equation}
where we put
\begin{equation}
\log\alpha_1^\prime=u_1\log\alpha_2-u_2\log\alpha_1,
\log\alpha_2^\prime=u_3\log\alpha_1+u_1\log\alpha_3.
\end{equation}

We recall that $\alpha_1=\si$, $\alpha_2=\eta_1$, $\alpha_3=\eta_2$ to obtain
\begin{equation}
\begin{split}
& \rho_1\abs{\log\alpha_i^\prime}+2\h(\alpha_i^\prime)
<u_1\left(\log m_i+\frac{\rho_1\pi}{2}\right)+\frac{u_{i+1}\rho_1 \pi}{2} \\
& <(1+\epsilon)\chi_1^{-2/3} C_1^{1/3}\left(\left(\log m_i+\frac{\rho_1\pi}{2}\right)\frac{\rho\pi}{2}
+\left(\log m_i+\frac{\rho\pi}{2}\right)\frac{\rho_1\pi}{2}\right) \\
& =(1+\epsilon)\chi_1^{-2/3} C_1^{1/3}((\rho+\rho_1)\log m_i+\rho\rho_1\pi)\pi/2.
\end{split}
\end{equation} 
So that we can take
\begin{equation}\label{eq36}
a_i^\prime=(1+\epsilon)\chi_1^{-2/3} C_1^{1/3}((\rho+\rho_1)\log m_i+\rho\rho_1\pi)\pi/2
\end{equation}
for $i=1, 2$ and therefore
\begin{equation}
\bpp=\frac{b_2}{a_2^\prime}+\frac{b_3}{a_1^\prime}<
\frac{2\chi_1^{2/3}}{\pi(\rho+\rho_1)C_1^{1/3}}
\left(\frac{e_1}{\log m_2}+\frac{e_2}{\log m_1}\right).
\end{equation}

Now we apply Lemma \ref{lm32} with $\rho_1$ taken from Table \ref{tbl2} and $\mu=0.61$ to obain
\begin{equation}
-\log\abs{\Lm}<A_1 \Omega \log^2 \bpp<A_1 \Omega \log^2(A_2^{\prpr} E),
\end{equation}
where $A_2^\prpr =2\chi_1^{2/3}/(\pi(\rho+\rho_1)C_1^{1/3})$.
Now we easily obtain \eqref{eq33} observing that $A_2^\prpr <A_2<A_2\sqrt{\log m_1}$.

Finally, if $u_1=0$, then $u_2 b_2+u_3 b_3=0$ and
therefore we can take $u_2=b_3^\prpr=e_2/d_2$ and $u_3=-b_2^\prime=-e_1/d_2$.
Hence,
\begin{equation}
\abs{b_2^\prime}\leq\frac{(R_{1, 1}+1)(R_{2, 1}+1)}{M-\max\{R_{1, 1}, R_{2, 1}\}}
<(1+\epsilon)\chi_1^{-2/3}C_1^{1/3}a_3
\end{equation}
and
\begin{equation}
\abs{b_3^\prpr}\leq\frac{(R_{1, 1}+1)(R_{3, 1}+1)}{M-\max\{R_{1, 1}, R_{3, 1}\}}
<(1+\epsilon)\chi_1^{-2/3}C_1^{1/3}a_2.
\end{equation}
Lemma \ref{lm22} immediately gives that $d_2=\gcd(e_1, e_2)\leq 4$ and
\begin{equation}
e_1<4(1+\epsilon)\chi_1^{-2/3}C_1^{1/3}\left(\frac{\rho\pi}{2} + \log m_2\right)
\end{equation}
and
\begin{equation}
e_2<4(1+\epsilon)\chi_1^{-2/3}C_1^{1/3}\left(\frac{\rho\pi}{2} + \log m_1\right).
\end{equation}
Hence, we obtain
\begin{equation}
E<4(1+\epsilon)\chi_1^{-2/3}C_1^{1/3}\left(2+\frac{\rho\pi}{2\log m_1}+\frac{\rho\pi}{2\log m_2}\right)<A_0.
\end{equation}
This completes the proof of the theorem.

For example, in the case $(\tilde m_1, \tilde m_2)=(53, m_1^*)$,
we have $K\geq 636133$, $L\geq 42$, and $R_i\geq 33242$ for $i=1, 2, 3$.
Hence, we obtain $R_1 R_2 R_3<18.58295K^2 L$, $g<0.2455157$, $C_3<31.908855$, and therefore
the right of \eqref{eq31b} is at most $46.93798$.
On the other hand, we have
\begin{equation}
C_0 C_1\left(1-\frac{1}{K}\right)\left(1-\frac{1}{L}-\frac{1}{3L^2}\right)>46.93857
\end{equation}
and we can confirm \eqref{eq31b}.

Now we can apply Lemma \ref{lm31}.
In the case (A), proceeding as above, we are led to
\begin{equation}
-\log\abs{\Lm}<5381.263\Omega\log^2 (1.985977E\sqrt{\log m_1}).
\end{equation}

In the case (B), we have from \eqref{eq34a} that
\begin{equation}
e_1\log m_1+e_2\log m_2<109.7095\Omega\log b_0<219.419\Omega\log(1.985977E\sqrt{\log m_1}).
\end{equation}

In the case (C),
applying Lemma \ref{lm32} with 
$\rho_1=6.2$ and $\mu=0.61$ and observing that
\begin{equation}
\bpp<\frac{2\chi_1^{2/3}}{C_1^{1/3}(\rho+\rho_1)\pi}\left(\frac{e_1}{\log m_1}+\frac{e_2}{\log m_2}\right)
<0.014249E,
\end{equation}
we obtain
\begin{equation}
-\log\abs{\Lm}<5381.263\Omega\log^2 \bpp
<5381.263\Omega\log^2 (0.014249E)
\end{equation}
or $E<A_0$.
These estimates give the theorem in the case $(\tilde m_1, \tilde m_2)=(53, m_1^*)$.
\end{proof}

\begin{table}
\caption{Other constants used in the proof of Theorem \ref{thm31}}
\begin{center}
\begin{small}
\begin{tabular}{| c | c | c | c | c | c | c | c | c |}
 \hline
$m_1$ & $\tilde m_2$ & $\chi_1$ & $C_0$ & $C_1$ & $\rho_1$ & $K_0$ & $L_0$ & $R_0$ \\
 \hline
$5$ & $2000000005$ & $0.46$ & $3.2$ & $15.499$ & $6.0$ & $549359$ & $42$ & $26154$ \\
$13$ & $769230773$ & $0.46$ & $3.3$ & $14.320$ & $6.1$ & $569527$ & $43$ & $28789$ \\
$17$ & $588235297$ & $0.50$ & $3.3$ & $11.653$ & $6.8$ & $583575$ & $40$ & $30786$ \\
$29$ & $344827589$ & $0.48$ & $3.3$ & $12.650$ & $6.5$ & $596786$ & $41$ & $31642$ \\
$37$ & $270270277$ & $0.53$ & $3.4$ & $9.571$ & $7.3$ & $602446$ & $39$ & $33973$ \\
$41$ & $243902441$ & $0.47$ & $3.3$ & $13.519$ & $6.2$ & $608540$ & $42$ & $32359$ \\
$5$ & $10^{16}$ & $0.48$ & $3.3$ & $13.487$ & $6.3$ & $822747$ & $42$ & $26520$ \\
$13$ & $10^{16}$ & $0.49$ & $3.3$ & $12.744$ & $6.4$ & $907924$ & $41$ & $29588$ \\
$17$ & $10^{16}$ & $0.54$ & $3.3$ & $10.207$ & $7.1$ & $943815$ & $38$ & $32235$ \\
$29$ & $10^{16}$ & $0.52$ & $3.3$ & $11.000$ & $6.9$ & $979765$ & $39$ & $32936$ \\
$37$ & $10^{16}$ & $0.55$ & $3.4$ & $9.013$ & $7.7$ & $984944$ & $38$ & $34781$ \\
$41$ & $10^{16}$ & $0.51$ & $3.3$ & $11.673$ & $6.6$ & $1008482$ & $40$ & $33538$ \\
$5$ & $10^{24}$ & $0.48$ & $3.3$ & $13.435$ & $6.4$ & $1163620$ & $42$ & $26532$ \\
$13$ & $10^{24}$ & $0.49$ & $3.3$ & $12.708$ & $6.5$ & $1282306$ & $41$ & $29594$ \\
$17$ & $10^{24}$ & $0.54$ & $3.4$ & $9.648$ & $7.3$ & $1290303$ & $39$ & $32099$ \\
$29$ & $10^{24}$ & $0.52$ & $3.4$ & $10.388$ & $7.1$ & $1340720$ & $40$ & $32759$ \\
$37$ & $10^{24}$ & $0.55$ & $3.5$ & $8.533$ & $7.9$ & $1345753$ & $39$ & $34693$ \\
$41$ & $10^{24}$ & $0.51$ & $3.3$ & $11.658$ & $6.7$ & $1421332$ & $40$ & $33543$ \\
$5$ & $e^{40000}$ & $0.51$ & $3.3$ & $11.847$ & $6.8$ & $739807372$ & $40$ & $27513$ \\
$13$ & $e^{40000}$ & $0.51$ & $3.4$ & $11.030$ & $7.0$ & $786646356$ & $41$ & $30127$ \\
$17$ & $e^{40000}$ & $0.59$ & $3.4$ & $8.245$ & $7.9$ & $817780793$ & $37$ & $34383$ \\
$29$ & $e^{40000}$ & $0.57$ & $3.4$ & $8.813$ & $7.6$ & $846771304$ & $38$ & $34859$ \\
$37$ & $e^{40000}$ & $0.64$ & $3.4$ & $7.019$ & $8.7$ & $888581986$ & $35$ & $39397$ \\
$41$ & $e^{40000}$ & $0.57$ & $3.3$ & $9.513$ & $7.3$ & $898313919$ & $37$ & $35840$ \\
 \hline
$\tilde m_1$ & $\tilde m_2$ & $\chi_1$ & $C_0$ & $C_1$ & $\rho_1$ & $K_0$ & $L_0$ & $R_0$ \\
 \hline
$53$ & \addstackgap[.5\dimexpr 2pt \relax]{$m_1^*$} &
$0.46$ & $3.2$ & $15.029$ & $6.2$ & $636133$ & $42$ & $33242$ \\
$53$ & $10^{24}$ & $0.49$ & $3.3$ & $12.550$ & $6.6$ & $1452862$ & $41$ & $33647$ \\
$53$ & $10^{32}$ & $0.51$ & $3.3$ & $11.688$ & $6.7$ & $1875132$ & $40$ & $34263$ \\
$53$ & $10^{34}$ & $0.51$ & $3.3$ & $11.683$ & $6.8$ & $1980530$ & $40$ & $34264$ \\
$53$ & $e^{2000}$ & $0.53$ & $3.4$ & $10.254$ & $7.1$ & $44504307$ & $40$ & $34762$ \\
$53$ & $e^{40000}$ & $0.53$ & $3.3$ & $10.813$ & $7.1$ & $913052819$ & $39$ & $34978$ \\
$10^{24}$ & $10^{24}$ & $0.64$ & $3.5$ & $6.286$ & $9.1$ & $6409655$ & $35$ & $176815$ \\
$5000$ & $e^{2000}$ & $0.59$ & $3.4$ & $8.211$ & $7.9$ & $59497348$ & $37$ & $49865$ \\
$5000$ & $e^{40000}$ & $0.59$ & $3.4$ & $8.203$ & $7.9$ & $1184137409$ & $37$ & $49870$ \\
$10^7$ & $e^{40000}$ & $0.63$ & $3.5$ & $6.784$ & $8.7$ & $1604506209$ & $36$ & $72224$ \\
$10^{14}$ & $e^{40000}$ & $0.69$ & $3.5$ & $5.745$ & $9.7$ & $2520730598$ & $34$ & $120077$ \\
 \hline
\end{tabular}
\label{tbl2}
\end{small}
\end{center}
\end{table}
\begin{table}
\caption{Other constants used in the proof of Theorem \ref{thm31}}
\begin{center}
\begin{small}
\begin{tabular}{| c | c | c | c | c | c | c |}
 \hline
$m_1$ & $\tilde m_2$ & $g_0$ & $C_{3, 0}$ & $X_0$ & $Y_0$ & $A_2^\prime$ \\
 \hline
$5$ & $2000000005$ & $0.2455091$ & $32.907152$ & $48.406238$ & $48.406463$ & $4.593026$ \\
$13$ & $769230773$ & $0.2454815$ & $31.827636$ & $46.147713$ & $46.148423$ & $4.481626$ \\
$17$ & $588235297$ & $0.2457133$ & $25.832053$ & $37.485114$ & $37.485451$ & $4.173153$ \\
$29$ & $344827589$ & $0.2456108$ & $28.067845$ & $40.717910$ & $40.718483$ & $4.291573$ \\
$37$ & $270270277$ & $0.2458429$ & $22.128118$ & $31.699167$ & $31.699820$ & $3.911263$ \\
$41$ & $243902441$ & $0.2455453$ & $30.022687$ & $43.541755$ & $43.541991$ & $4.392677$ \\
$5$ & $10^{16}$ & $0.2455813$ & $29.951582$ & $43.438633$ & $43.438944$ & $4.393314$ \\
$13$ & $10^{16}$ & $0.2456439$ & $28.276664$ & $41.020708$ & $41.021078$ & $4.307492$ \\
$17$ & $10^{16}$ & $0.2459050$ & $22.581641$ & $32.788675$ & $32.788892$ & $3.991684$ \\
$29$ & $10^{16}$ & $0.2458083$ & $24.360773$ & $35.360809$ & $35.361239$ & $4.094194$ \\
$37$ & $10^{16}$ & $0.2459353$ & $20.816751$ & $29.829780$ & $29.830669$ & $3.833814$ \\
$41$ & $10^{16}$ & $0.2457464$ & $25.876539$ & $37.549576$ & $37.549815$ & $4.179917$ \\
$5$ & $10^{24}$ & $0.2455822$ & $29.835631$ & $43.270668$ & $43.271477$ & $4.387128$ \\
$13$ & $10^{24}$ & $0.2456445$ & $28.197090$ & $40.905122$ & $40.905213$ & $4.303060$ \\
$17$ & $10^{24}$ & $0.2458731$ & $22.305982$ & $31.954006$ & $31.954878$ & $3.926510$ \\
$29$ & $10^{24}$ & $0.2457767$ & $24.040295$ & $34.428320$ & $34.428836$ & $4.026041$ \\
$37$ & $10^{24}$ & $0.2459042$ & $20.559644$ & $29.092665$ & $29.093151$ & $3.772791$ \\
$41$ & $10^{24}$ & $0.2457467$ & $25.843209$ & $37.501236$ & $37.501573$ & $4.177959$ \\
$5$ & $e^{40000}$ & $0.2457432$ & $26.261907$ & $38.108908$ & $38.109577$ & $4.202487$ \\
$13$ & $e^{40000}$ & $0.2457150$ & $25.549364$ & $36.579365$ & $36.579880$ & $4.111220$ \\
$17$ & $e^{40000}$ & $0.2460916$ & $19.023386$ & $27.268387$ & $27.268525$ & $3.725578$ \\
$29$ & $e^{40000}$ & $0.2460016$ & $20.355247$ & $29.168248$ & $29.168751$ & $3.810001$ \\
$37$ & $e^{40000}$ & $0.2463000$ & $16.156984$ & $23.175985$ & $23.176260$ & $3.529766$ \\
$41$ & $e^{40000}$ & $0.2460292$ & $21.023538$ & $30.536539$ & $30.536799$ & $3.901703$ \\
 \hline
$\tilde m_1$ & $\tilde m_2$ & $g_0$ & $C_{3, 0}$ & $X_0$ & $Y_0$ & $A_2^\prime$ \\
 \hline
$53$ & \addstackgap[.5\dimexpr 2pt \relax]{$m_1^*$} &
$0.2455157$ & $31.908855$ & $46.937927$ & $46.938571$ & $4.541758$ \\
$53$ & $10^{24}$ & $0.2456472$ & $27.846484$ & $40.396512$ & $40.396637$ & $4.283496$ \\
$53$ & $10^{32}$ & $0.2457461$ & $25.909877$ & $37.597898$ & $37.598084$ & $4.181873$ \\
$53$ & $10^{34}$ & $0.2457462$ & $25.898758$ & $37.581778$ & $37.582001$ & $4.181221$ \\
$53$ & $e^{2000}$ & $0.2458120$ & $23.730227$ & $33.984228$ & $33.984745$ & $4.010986$ \\
$53$ & $e^{40000}$ & $0.2458444$ & $23.947111$ & $34.760112$ & $34.760133$ & $4.072742$ \\
$10^{24}$ & $10^{24}$ & $0.2462967$ & $15.080266$ & $21.366270$ & $21.366409$ & $3.403729$ \\
$5000$ & $e^{2000}$ & $0.2460924$ & $18.944866$ & $27.155867$ & $27.156077$ & $3.719988$ \\
$5000$ & $e^{40000}$ & $0.2460926$ & $18.926542$ & $27.129543$ & $27.129619$ & $3.718673$ \\
$10^7$ & $e^{40000}$ & $0.2462366$ & $16.294300$ & $23.078301$ & $23.078337$ & $3.496385$ \\
$10^{14}$ & $e^{40000}$ & $0.2464633$ & $13.764710$ & $19.509711$ & $19.510304$ & $3.309117$ \\
\hline
\end{tabular}
\label{tbl3}
\end{small}
\end{center}
\end{table}

Furthermore, we need a lower bound for linear forms in one logarithm and $\pi \si/2$.
\begin{lemma}\label{lm33}
Let
\begin{equation}
\Lm_1=b_2\log \alpha-\frac{b_1\pi \si}{2},
\end{equation}
where $b_1, b_2$ are positive integers,
$\alpha$ is an complex algebraic number of absolute value one but not a root of unity
and $\log\alpha$ takes its value such that $\abs{\Imf \log\alpha}<\pi/2$.
Put
\begin{equation}
\begin{split}
\bp= & \frac{b_1}{9\pi+2D\h(\alpha)}+\frac{b_2}{9\pi}, \\
D= & [\Q(\alpha): \Q]/2, \\
a= & 9\pi+2D\h(\alpha), \\
H= & \max\{17, D, D(\log \bp+2.96)+0.01\}.
\end{split}
\end{equation}
Then,
\begin{equation}
\log\abs{\Lm_1}>-(2.76701a+0.12945)H^2.
\end{equation}
\end{lemma}

\begin{proof}
This follows from Theorem 1 of \cite{Lau} with $D=1$ by setting
$\mu=0.59$, $\rho_1=18$, $\alpha_1=i$, $\alpha_2=\alpha$,
$a_1=\rho\pi/2$, $a_2=2D\h(\alpha)+(\rho\pi/2)$, $L=\floor{H+\sq{H^2+(1/4)}+(1/2)}$,
$k=(v_1(L)+\sq{v_1(L)^2+4v_0(L)v_2(L)})^2/(2v_2(L))^2$ to be the positive real number such that
$\sq{k}$ satisfies the quadratic equation
$v_2(L) k-v_1(L)\sq{k}-v_0(L)=0$,
where
$$
v_0(x)=\frac{1}{4a_1}+\frac{4}{3a_2}+\frac{x}{12a_1},
v_1(x)=\frac{x}{3},
v_2(x)=\lm (x-H),
$$
and $K=1+\floor{kLa_1 a_2}$, $R_1=4$, $S_1=\floor{(L+3)/4}$, 
$R_2=1+\floor{\sq{(K-1)La_2/a_1}}$, $S_2=1+\floor{\sq{(K-1)La_1/a_2}}$.

Then, we can confirm that $0.2437<\sq{k}<0.2791$, $H>7.5$, $15\leq L<0.92h$, and $K>712$
and then the condition (2) of Theorem 1 of \cite{Lau}.
Moreover, we can also confirm that $R_2-1<2a_2 H^2<b_1 $ and $S_2-1<2a_1 H^2<b_2$.
Hence the condition (1) of Theorem 1 of \cite{Lau} also holds.
Then, confirming that $\sqrt{k}L<0.239537H$, Theorem 1 of \cite{Lau} yields that
$\log\abs{\Lm_1\max\{LS e^{LS\abs{\Lm_1}/(2b_2)}/(2b_2), LR e^{LR\abs{\Lm_1}/(2b_1)}/(2b_1)\}}>-(2.7666a+0.10984)H^2$
and then $\log\abs{\Lm_1}>-(2.76701a+0.12945)aH^2$.
For the detail of the proof, see \cite{Ymd2}.
\end{proof}

\section{Upper bounds for exponents}

In this section, we shall prove upper bounds for exponents in \eqref{eq12} or, equivalently, \eqref{eq13}.

\begin{theorem}\label{thm41}
Let $\eta_1$ and $\eta_2$ be Gaussian integers with $-\pi/2<\arg \eta_i/\bar\eta_i<\pi/2$
such that $m_1=\eta_1\bar\eta_1$ and $m_2=\eta_2\bar\eta_2$ satisfy the condition of Theorem \ref{thm31}.
If $x, e_1, e_2$ are nonnegative integers such that $x^2+1=m_1^{e_1} m_2^{e_2}$
or $x^2+1=2m_1^{e_1} m_2^{e_2}$, then,
\begin{equation}\label{eq41}
e_1\log m_1+e_2\log m_2\leq A_4\Omega\log^2 (A_5 \log m_1),
\end{equation}
where $A_5=4A_2^2(A_1+0.6932/\log^2(10000A_2\sqrt{\log \tilde m_1}))^2$,
$A_4$ are constants given in Table \ref{tbl4},
and $\Omega=(\log m_1)(\log m_2)$ as in the previous section.
Moreover, if $m_1 m_2<10^{10}$, then we have
\begin{equation}\label{eq41b}
e_1\log m_1+e_2\log m_2<1.7628 m_1 m_2.
\end{equation}
\end{theorem}

\begin{table}
\caption{Constants in \eqref{eq41}}
\begin{center}
\begin{small}
\begin{tabular}{| c | c | c |}
 \hline
$m_1$ & $\tilde m_2$ & $A_4$ \\
 \hline
$5$ & $2000000005$ & $7959.719$ \\
$13$ & $769230773$ & $6815.554$ \\
$17$ & $588235297$ & $3981.632$ \\
$29$ & $344827589$ & $4533.458$ \\
$37$ & $270270277$ & $3127.262$ \\
$41$ & $243902441$ & $5922.421$ \\
$5$ & $10^{16}$ & $7192.899$ \\
$13$ & $10^{16}$ & $6087.270$ \\
$17$ & $10^{16}$ & $3510.530$ \\
$29$ & $10^{16}$ & $3979.106$ \\
$37$ & $10^{16}$ & $2698.486$ \\
$41$ & $10^{16}$ & $5206.488$ \\
$5$ & $10^{24}$ & $6811.925$ \\
$13$ & $10^{24}$ & $5757.606$ \\
$17$ & $10^{24}$ & $3305.365$ \\
$29$ & $10^{24}$ & $3754.144$ \\
$37$ & $10^{24}$ & $2533.584$ \\
$41$ & $10^{24}$ & $4914.868$ \\
$5$ & $e^{40000}$ & $6033.728$ \\
$13$ & $e^{40000}$ & $5083.203$ \\
$17$ & $e^{40000}$ & $2875.896$ \\
$29$ & $e^{40000}$ & $3281.774$ \\
$37$ & $e^{40000}$ & $2176.067$ \\
$41$ & $e^{40000}$ & $4350.349$ \\
 \hline
\end{tabular}
\begin{tabular}{| c | c | c |}
 \hline
$\tilde m_1$ & $\tilde m_2$ & $A_4$ \\
 \hline
$53$ & \addstackgap[.5\dimexpr 2pt \relax]{$m_1^*$} & $6240.476$ \\
$53$ & $10^{24}$ & $5167.476$ \\
$53$ & $10^{32}$ & $5024.642$ \\
$53$ & $10^{34}$ & $4997.394$ \\
$53$ & $e^{2000}$ & $4579.094$ \\
$53$ & $e^{40000}$ & $4555.827$ \\
$10^{24}$ & $10^{24}$ & $1752.029$ \\
$5000$ & $e^{2000}$ & $2889.083$ \\
$5000$ & $e^{40000}$ & $2876.189$ \\
$10^7$ & $e^{40000}$ & $2131.592$ \\
$10^{14}$ & $e^{40000}$ & $1675.249$ \\
 \hline
\end{tabular}
\label{tbl4}
\end{small}
\end{center}
\end{table}

\begin{proof}
Assume that $m_1$ and $m_2$ satisfy the condition of Theorem \ref{thm31}
and put $\tilde m_1$, $\tilde m_2$, $\Omega$, and $E$ as in Theorem \ref{thm31}.
We may assume that $e_1 e_2\neq 0$ since $e_i=0$ implies that $e_{3-i}=1, 2$ or $4$
by Lemma \ref{lm22}, so that we can put $\alpha=(e_2\log m_2)/(e_1\log m_1)$.
Moreover, we may assume that $E\geq 10000$ since otherwise we have
$e_1\log m_1+e_2\log m_2<10000\Omega<A_4\Omega \log^2 A_5$.

We can decompose $m_i=\eta_i \bar\eta_i$ in a way such that $-\pi/4<\arg \eta^\prime<\pi/4$.
We put $\xi_i=\eta_i/\bar\eta_i$ and write $\theta_i=\abs{\arg\xi_i}=\abs{\log \xi_i}$,
so that $\theta_i<\pi/2$.

Now $\Lm=\log[(x+\si)/(x-\si)]$ can be represented as a linear form of three logarithms
\begin{equation}
\Lm=\pm e_1\log \xi_1\pm e_2\log \xi_2\pm \frac{e_3\pi \si}{2}
\end{equation}
for an appropriate integer $e_3\geq 0$.
Moreover, we can easily see that
\begin{equation}
\log \abs{\Lm}<-\log (x/2)<-\frac{e_1\log m_1 + e_2\log m_2}{2}+10^{-9}+\log 2=-\frac{E\Omega}{2}+10^{-9}+\log 2.
\end{equation}
Since $m_1$ and $m_2$ satisfy the condition of Theorem \ref{thm31},
we use this theorem to obtain either (a) \eqref{eq33}, (b) \eqref{eq34},
or (c) $(A_2 E)^2\log m_1\leq B_0$.

If \eqref{eq33} holds, then
\begin{equation}
\frac{E\Omega}{2}<0.6932+A_1 \log^2 (A_2 E\sqrt{\log m_1}) \Omega.
\end{equation}
Hence, we have
\begin{equation}
A_2 E\sqrt{\log m_1}<2A_2 \sqrt{\log m_1} (0.6932+A_1\log^2 (A_2 E\sqrt{\log m_1}))
\end{equation}
and therefore $A_2 E\sqrt{\log m_1}<\delta \log^2 (A_2 E\sqrt{\log m_1})$
with $\delta=\sqrt{A_5\log m_1}$.

Let $c_5$ be the unique positive number satisfying
\begin{equation}
c_5=\left(1+\frac{2\log\log (\sqrt{A_5 \log \tilde m_1})+\log c_5}
{\log (\sqrt{A_5 \log \tilde m_1})}\right)^2.
\end{equation}
Since $\delta\geq A_5\sqrt{\log \tilde m_1}$, we have
\begin{equation}
c_5\geq \left(1+\frac{2\log\log\delta+\log c_5}{\log\delta}\right)^2.
\end{equation}
Now, if $X>c_5\delta\log^2\delta$, then
\begin{equation}
\frac{X}{\log^2 X}>\frac{c_5\delta\log^2\delta}{(\log\delta+2\log\log\delta+\log c_5)^2}
=\frac{c_5\delta}{\left(1+\frac{2\log\log\delta +\log c_5}{\log\delta}\right)^2}\geq \delta.
\end{equation}
Hence, we must have $A_2 E\sqrt{\log m_1}<c_5 \delta\log^2 \delta$ and
\begin{equation}
E<\frac{c_5}{2} (A_1+0.6932/\log^2(10000A_2\sqrt{\log \tilde m_1}))\log^2(A_5 \log m_1).
\end{equation}
Since $A_4\geq c_5 (A_1+0.6932/\log^2(10000A_2\sqrt{\log \tilde m_1}))/2$, \eqref{eq41} holds.

If \eqref{eq34} holds, then we have
\begin{equation}
A_2 E\sqrt{\log m_1}<A_2 A_3 \log (A_2 E\sqrt{\log m_1})\sqrt{\log m_1}
\end{equation}
and, proceeding like above, we have $E<A_4 \log^2(A_5 \log m_1)$, which immediately yields \eqref{eq41}.

If $(A_2 E)^2\log m_1\leq B_0$, then we immediately have $A_2 E<A_2 E\sqrt{\log m_1}<\sqrt{B_0}$ and
$E<\sqrt{B_0}/A_2$.
But we can confirm that $\sqrt{B_0}/A_2<A_4\log^2 A_5$ in any case.
Hence, we obtain \eqref{eq41}.

Finally, we assume that $m_1 m_2<10^{10}$.
We observe that there exist two squarefree integers $D_1$ and $D_2$ such that $D_1 D_2$ divides $2m_1 m_2$ and
$x^2-D_1 y^2=-1$
for some integer $y$ all of whose prime factors divide $D_2$.
We have $x+y\sqrt{D_1}=\epsilon^n$ for some integer $n$, where $\epsilon$ is the fundamental unit in $\Q(D_1)$
Now it follows from Proposition 3.4.5 of \cite[p. 138]{HC1} and Proposition 10.3.16 of \cite[p. 200]{HC2}
with $D_1$ in place of $m$ and $f$ in the latter proposition that
\begin{equation}
\log\epsilon<D_1^{1/2}((\log D_1)/2+\log \log D_1+2.8)<0.4407 D_1.
\end{equation}
By Carmichael's theorem \cite{Car} (for a simple proof, see \cite{Ybt}),
we see that $n$ divides $p-(D_1/p)$ for some prime factor $p$ of $y$ (and therefore $D_2$) whenever $n>12$.
Hence, we have $n\leq D_2+1$ to obtain
\begin{equation}
\log(x+y\sqrt{D_1})=n\log\epsilon<0.4407D_1(D_2+1)<0.8814m_1 m_2,
\end{equation}
which yields that
\begin{equation}
e_1\log m_1+e_2\log m_2\leq \log(x^2+1)=\log(D_1 y^2)<1.7628 m_1 m_2,
\end{equation}
which completes the proof of the lemma.
\end{proof}

\section{Proof of the Theorem}

Let $x_1, x_2, x_3, y_1, y_2, y_3$ and $r$ be integers with $x_1, x_2, x_3>1, r\neq 0$ and
$\gcd(y_1, y_2, y_3)=1$ satisfying \eqref{eq11}
and $m_1, m_2, s_i, k_i, \ell_i (i=1, 2, 3)$ be corresponding integers.
Furthermore, let $\eta_1, \eta_2$ be Gaussian integers satisfying \eqref{eq12} and \eqref{eq13}.
We write $F=\max k_i$, $G=\max\ell_i$, and, as in preceding sections, $\Omega=(\log m_1)(\log m_2)$.
Moreover, we write the quantity $\max\{A, B\}$ as $A\lor B$ for brevity.

We note that, since we have assumed that $\{x_1, x_2, x_3\}\neq \{2, 3, 7\}$,
\eqref{eq11} is nondegenerate by Lemma \ref{lm21} and therefore
we have $y_1\mid (\pm k_2 \ell_3\pm k_3 \ell_2), y_2\mid (\pm k_3 \ell_1\pm k_1 \ell_3)$ and
$y_3\mid (\pm k_1 \ell_2\pm k_2 \ell_1)$ with appropriate choices of signs.
Hence, $\abs{y_i}\leq 2FG$ for $i=1, 2, 3$.

Now we assume that $m_1 m_2\geq 10^{10}$.
Then, the former part of Theorem \ref{thm41} yields that, for each $i$, 
\begin{equation}\label{eq51a}
k_i\leq A_4(\log m_2)\log^2 (A_5 \log m_1),
\end{equation}
where we observe that $A_4^2\log^4 A_5>e^{20}>B_0/A_2^2$ in any case, and, similarly,
\begin{equation}\label{eq51b}
\ell_i\leq A_4(\log m_1)\log^2 (A_5 \log m_1).
\end{equation}
Hence, we have
\begin{equation}\label{eq51}
FG\leq A_4^2\Omega\log^4 (A_5 \log m_1).
\end{equation}

We have two cases: I. $m_2$ divides $x_i^2+1$ for $i=1, 2, 3$,
and II. $x_i^2+1=m_1^{k_i}$ or $2m_1^{k_i}$ for some $i$.

\textbf{Case I.}
In this case, we have two possibilities (a) $x_i^2+1=m_2^{\ell_i}$ or $2m_2^{2\ell_i}$ for some $i$
and (b) $m_1 m_2$ divides $x_i^2+1$ for each $i=1, 2, 3$.

In the case (a), we arrange indices so that $x_1^2+1=m_2^{\ell_1}$ or $2m_2^{\ell_1}$.
By Lemma \ref{lm22}, we have $\ell_1=1, 2$, or $(x_1, m_2, \ell_1)=(239, 13, 4)$.
In the latter case, we can replace $(x_1, \ell_1)$ by $(5, 1)$.
Hence, we may assume that $\ell_1=1$ or $2$.

Since $x_1^2+1\geq m_2$, we have $x_1\geq \sqrt{m_2-1}$.
Since $x_2^2+1\equiv x_3^2+1\equiv 0\Mod{m_2}$, St{\o}rmer's criterion implies that
$x_2>m_2/2$ and $x_3>m_2$.

Since $\abs{y_2}\leq \abs{k_3\ell_1}+\abs{k_1\ell_3}\leq 2\abs{k_3}$ and
$\abs{y_3}\leq \abs{k_2\ell_1}+\abs{k_1\ell_2}\leq 2\abs{k_2}$, we must have
\begin{equation}
\frac{\abs{y_1}}{\sqrt{m_2-1}}+\frac{2(2\abs{k_2}+\abs{k_3})}{m_2}>\frac{\pi}{4}
\end{equation}
and therefore
\begin{equation}
\frac{2FG}{\sqrt{m_2-1}}+\frac{6F}{m_2}>\frac{\pi}{4}.
\end{equation}
The assumption that $m_1 m_2>10^{10}$ implies that $m_2>10^5$ and therefore
$(2.000011FG+0.018974F)/\sqrt{m_2}>\pi/4$.
Hence, we have $m_2<(2.5465FG+0.02416F)^2$.

In the case (b), we can arrange indices so that $x_3>x_2>x_1$.
Since $x_1^2+1\geq m_1 m_2\geq 5m_2$, we have $x_1\geq \sqrt{5m_2-1}$.
Since $x_2^2+1\equiv x_3^2+1\equiv 0\Mod{m_2}$, St{\o}rmer's criterion implies that
$x_2>m_2/2$ and $x_3>m_2$.

It immediately follows from \eqref{eq11} with $r\neq 0$ that
\begin{equation}
\frac{\abs{y_1}}{\sqrt{5m_2-1}}+\frac{2\abs{y_2}+\abs{y_3}}{m_2}>\frac{\pi}{4}.
\end{equation}
Since $\abs{y_i}\leq 2FG$ for $i=1, 2, 3$, we have
\begin{equation}
\frac{2FG}{(\sqrt{5-10^{-6}})\sqrt{m_2}}>\frac{2FG}{\sqrt{5m_2-1}}+\frac{6FG}{m_2}>\frac{\pi}{4}.
\end{equation}
and therefore $m_2<1.306(FG)^2$.
Hence, in both subcases, we have $m_2<(2.5465FG+0.02416F)^2$.

On the other hand, we have
\begin{equation}
FG<A_4^2(\log^2 m_2)\log^4 (A_5 \log m_2)
\end{equation}
and
\begin{equation}
F<A_4(\log m_2)\log^2 (A_5\log m_2).
\end{equation}

If $m_1\geq 53$ and $m_2>1.342\times 10^{34}$, then, we can take $A_4=4997.394$ from Table \ref{tbl4}.
Moreover, we have $A_5<168193563$.
Hence, we obtain
\begin{equation}
m_2<4.0444864\times 10^{15}\log^8 (168193563\log m_2) \log^4 m_2.
\end{equation}
However, this cannot hold for $m_2>1.342\times 10^{34}$.
Hence, we must have $m_1<m_2<1.342\times 10^{34}$ when $m_1\geq 53$.

If $m_2\geq 10^{34}$, then, applying Theorem \ref{thm41} with $(\tilde m_1, \tilde m_2)=(53, 10^{34})$ again,
we have $\abs{y_i}<9.098\times 10^{16}$ and $\log x_i<8380114833$.
If $10^{32}\leq m_2<10^{34}$, then, applying Theorem \ref{thm41} with $(\tilde m_1, \tilde m_2)=(53, 10^{32})$,
we have $\abs{y_i}<9.152\times 10^{16}$ and $\log x_i<8373165769$.
If $10^{24}\leq m_2<10^{32}$, then, applying Theorem \ref{thm41} with $(\tilde m_1, \tilde m_2)=(53, 10^{24})$,
we have $\abs{y_i}<8.752\times 10^{16}$ and $\log x_i<7706422610$.
If $m_2<10^{24}$, then, applying Theorem \ref{thm41} again with $\tilde m_1=53$ and $\tilde m_2=m_1^*$, we have $\abs{y_i}<7.799\times 10^{16}$ and $\log x_i<5456058110$. 
Hence, in any case, we have $y_i<9.152\times 10^{16}$ and $x_i<\exp(8380114833)$.

\begin{table}
\caption{Upper bounds for $\log x_i$ and $y_i ~ (i=1, 2, 3)$ (Case I)}
\begin{center}
\begin{small}
\begin{tabular}{| c | c | c | c | c | c |}
\hline
$m_1$ & $m_2<$ & $\log x_i<$ & $y_i<$ \\
\hline
$\geq 53$ & \addstackgap[.5\dimexpr 2pt \relax]{$1.342\times 10^{34}$} &
$8380114833$ & $9.152\times 10^{16}$ \\
$41$ & $7.462\times 10^{30}$ &
$266070791$ & $2.146\times 10^{15}$ \\
$37$ & $2.157\times 10^{29}$ &
$108769740$ & $3.647\times 10^{14}$ \\
$29$ & $1.448\times 10^{30}$ &
$185150405$ & $9.451\times 10^{14}$ \\
$17$ & $4.753\times 10^{29}$ &
$154361476$ & $5.415\times 10^{14}$ \\
$13$ & $6.816\times 10^{30}$ &
$324293344$ & $2.051\times 10^{15}$ \\
$5$ & $5.074\times 10^{30}$ &
$396883453$ & $1.769\times 10^{15}$ \\
\hline
\end{tabular}
\label{tbl5}
\end{small}
\end{center}
\end{table}

Similarly, if $m_1=41, 37, 29, 17, 13$, and $5$, then $m_2<7.461\times 10^{30}$,
$2.146\times 10^{29}$, $1.448\times 10^{30}$, $4.753\times 10^{29}$, $6.813\times 10^{30}$,
and $5.072\times 10^{30}$, respectively.

Now, applying Theorem \ref{thm41} again with $\tilde m_2=m_1^*$, we have
$\abs{y_i}\leq 2FG<2.146\times 10^{15}$
and
$\log x_i<(\log 2+k_i\log m_1+\ell_i\log m_2)/2<266065055$,
that is, $x_i<\exp(6263314097)$ for $m_1\leq 41$ (see Table \ref{tbl5}).
This shows the Theorem in Case I.

\textbf{Case II.}
We can arrange the indices so as to have $x_1^2+1=2^{v_1} m_1^{k_1}$ with $v_1=0$ or $1$ and $x_3>x_2$.
By Lemma \ref{lm22}, we see that $k_1=1$, $2$, or $4$ and $k_1=4$ implies $x_1=239$ and $m_1=13$ like above.
If $\ell_i=0$ for another index $i>1$, then Lemma \ref{lm22} implies that
$\{x_1, x_2, x_3\}=\{2, 3, 7\}$.
Thus, we must have $\ell_2, \ell_3>0$.

We begin by showing that $m_1<M_1$ with $M_1=2.531\times 10^{24}$.
Assume that $m_2\geq m_1\geq M_1$.
Now, Theorem \ref{thm41} allows us to take $A_4=1752.029$,
$A_5=10545512.8\cdots$, and $\rho=7.292$.

We shall settle the case $k_1\leq 2$, so that we leave the case $k_1=4$ for later.
Hence, we obtain
\begin{equation}\label{eq53a}
\begin{split}
FG<A_4^2\Omega\log^4 (10545513\log m_1).
\end{split}
\end{equation}
Similarly, we obtain
\begin{equation}\label{eq53b}
F<A_4(\log m_2)\log^2 (10545513 \log m_1)
\end{equation}
and
\begin{equation}\label{eq53c}
G<A_4(\log m_1)\log^2 (10545513 \log m_1).
\end{equation}
Hence, we obtain
\begin{equation}\label{eq53d}
\left(\frac{G\log m_2}{\log m_1}\lor F\right)<A_4 (\log m_2) \log^2 (10545513 \log m_1).
\end{equation}

Now, if $m_1\leq 7.2\times 10^8 G^2$, then, \eqref{eq53c} gives
$m_1<7.2\times 10^8 A_4^2(\log^2 m_1)\log^4 (10545513\log m_1)$ and therefore we must have
$m_1<1.2\times 10^{24}<M_1$.
Thus, we may assume that $m_2>m_1>7.2\times 10^8 G^2$.

We clearly have
\begin{equation}\label{eq54}
y_1 \arctan \frac{1}{x_1}\pm \ell_3 k_1 \arctan \frac{1}{x_2} \pm \ell_2 k_1 \arctan \frac{1}{x_3}=\frac{r\pi}{4}
\end{equation}
with $k$ some integer and $y_1=\ell_2 k_3\pm \ell_3 k_2$.
We can easily see that $x_2\geq \sqrt{m_1^{k_2} m_2^{\ell_2}-1}\geq \sqrt{m_2-1}$.
Moreover, from St{\o}rmer's criterion, we observe that $x_3>m_2/2$.
Since we have just assumed that $m_2\geq 7.2\times 10^8 G^2$, we obtain
\begin{equation}\label{eq55}
\abs{\frac{r\pi}{4}-y_1 \arctan \frac{1}{x_1}}<
\frac{k_1 \ell_3}{x_2-1}+\frac{k_1 \ell_2}{x_3-1}
<\frac{(2+10^{-5})G}{\sqrt{m_2}}<7.5\times 10^{-5}.
\end{equation}

Like above, we see that $\arctan(1/x_1)<1/(x_1-1)<1/((1-10^{-9})\sqrt{m_1})$ and therefore
\begin{equation}
\frac{\abs{y_1}}{\sqrt{m_1}}>(1-10^{-9})\left(\frac{\pi}{4}-7.5\times 10^{-5}\right)>0.785323.
\end{equation}
Since $k_2\log m_1+\ell_2 \log m_2\leq\log(x_2^2+1)$, putting $\alpha=\log(x_2^2+1)/\log m_2\geq 1$,
we obtain
\begin{equation}\label{eq56}
\begin{split}
\abs{y_1}\leq & ~ k_2 G+\ell_2 F\leq (\log (x_2^2+1))\left(\frac{G}{\log m_1}\lor \frac{F}{\log m_2}\right) \\
= & ~ \alpha\left(\frac{G\log m_2}{\log m_1}\lor F\right)
<A_4 (\alpha \log m_2)\log^2 (10545513 \log m_1)
\end{split}
\end{equation}
and therefore
\begin{equation}\label{eq57}
A_4 (\alpha \log m_2) \log^2 (10545513 \log m_1)>0.785319\sqrt{m_1}.
\end{equation}

Let
\begin{equation}
\Lm_1=y_1\log\frac{x_1+\si}{x_1-\si}-\frac{r\pi \si}{2}.
\end{equation}
Then \eqref{eq55} immediately gives that
\begin{equation}\label{eq55b}
\abs{\Lm_1}<\frac{(4+10^{-9})G}{m_1^{k_2/2} m_2^{\ell_2/2}},
\end{equation}
while Lemma \ref{lm33} gives that
\begin{equation}\label{eq58}
\log \abs{\Lm_1}>-(2.76701a+0.12945)aH^2,
\end{equation}
where
\begin{equation}
\begin{split}
a= & 9\pi+\log m_1, \\
H= & \max\left\{17, 2.97+\log \left(\frac{\abs{r}}{9\pi+\log m_1}+\frac{\abs{y_1}}{9\pi}\right)\right\}.
\end{split}
\end{equation}

It immediately follows from \eqref{eq55} and \eqref{eq56} that
\begin{equation}
\abs{\frac{r\pi}{4}}-7.5\times 10^{-5}\leq \frac{\abs{y_1}}{x_1-1}
<\frac{A_4(\alpha\log m_2) \log^2 (10545513 \log m_1)}{m_1^{1/2}}.
\end{equation}
Hence, we obtain
\begin{equation}
\frac{\abs{r}}{9\pi+\log m_1}+\frac{\abs{y_1}}{9\pi}<A_4\psi_3(m_1)(\alpha \log m_2)\log^2(10545513\log m_1),
\end{equation}
where
\begin{equation}
\begin{split}
\psi_3(m_1)= & ~ \frac{4}{(9\pi+\log m_1)\pi m_1^{1/2}}+\frac{1}{9\pi} \\
& ~ +\frac{3.16\times 10^{-4}}{\pi(9\pi+\log m_1)A_4(\alpha \log m_2)\log^2(10545513\log m_1)}.
\end{split}
\end{equation}
Since $m_1\geq M_1$, we have $2.97+\log\psi_3(m_1)\leq 2.97+\log\psi_3(M_1)<-0.377494$.
Hence, if $H>17$, then
\begin{equation}
\begin{split}
H< & ~ 2.97+\log\psi_3(m_1)+\log(A_4\alpha(\log m_2)\log^2 (10545513\log m_1)) \\
< & ~ -0.377494+\log(A_4\alpha (\log m_2)\log^2 (10545513\log m_1)).
\end{split}
\end{equation}
We observe that the right hand side of this inequality is at least \\
$-0.377494+\log(A_4\log M_1\log^2(10545513\log M_1))>17$
to obtain
\begin{equation}
H<-0.377494+\log(A_4\alpha (\log m_2)\log^2 (10545513\log m_1))
\end{equation}
whether $H>17$ or not.
With the aid of these inequalities, we deduce from \eqref{eq55b} and \eqref{eq58} that
\begin{equation}\label{eq59}
\begin{split}
& \alpha\log m_2<2(\log ((4+10^{-9})G)+(2.76701a+0.12945)H^2) \\
< & ~ 2 \log ((4+10^{-9})A_4\log^2 (10545513\log m_1)\log m_1) \\
& ~ +(5.53402\log m_1+156.73)(-0.377494+\log(A_4 \alpha(\log m_2) \log^2(10545513\log m_1)))^2.
\end{split}
\end{equation}

Recalling that $m_1, m_2\geq M_1$ and $\alpha\geq 1$, we must have
\begin{equation}\label{eq5a}
\alpha\log m_2<150.5379(\log m_1)\log^2 (\alpha\log m_2)
\end{equation}
and therefore
\begin{equation}\label{eq5b}
\alpha\log m_2<380.3347(\log m_1)(\log(150.5379\log m_1))^2<1916.227(\log m_1)(\log\log m_1)^2.
\end{equation}
This implies that
\begin{equation}\label{eq5c}
A_4 (\alpha \log m_2) \log^2 (10545513\alpha \log m_1)<0.785323\sqrt{m_1}
\end{equation}
when $m_1\geq M_1$, which is incompatible with \eqref{eq57}.
Thus, we are led to the estimate $m_1<M_1$.

Now we show that $m_2<\exp(294622)$.
Assume that $m_2\geq \exp(294622)$.
If $m_1\geq M_1^\prime:=9.134\times 10^{22}$, then, since $m_2>e^{40000}$,
we iterate the argument that we have just done with $(\tilde m_1, \tilde m_2)=(10^{14}, e^{40000})$
to obtain \eqref{eq5a} and \eqref{eq5b}
with $42.16707$, $116.41714$, and $439.5026$ in place of $150.5379$, $380.3347$, and $1916.227$.
Eventually, we obtain
\begin{equation}\label{eq5c0}
A_4 (\alpha \log m_2) \log^2 (6665783\alpha \log m_1)<0.785376\sqrt{m_1}
\end{equation}
with $A_4=1675.249$ instead of \eqref{eq5c}.
On the other hand, if $m_1\leq 5\times 10^{14} G^2$, then, \eqref{eq53c} gives
$G<A_4(\log^2 (5\times 10^{14} G^2))\log^4 (6665783\log (2.56\times 10^{10} G^2))$
and therefore $G<3\times 10^{12}$, which implies that $m_2>e^{40000}>5\times 10^{14} G^2$.
Hence, in any case, we must have $m_2>5\times 10^{14} G^2$ to obtain \eqref{eq55} and then \eqref{eq57}
with $10^{-6}$ in place of $7.5\times 10^{-5}$ and $0.785376$ in place of $0.785323$, contrary to \eqref{eq5c0}.
Hence, we must have $m_1<M_1^\prime$ when $m_2\geq e^{40000}$.

However, \eqref{eq59} still holds with $2.97+\log\psi_3(m_1)$ in place of $-0.377494$
and $3.16\times 10^{-4}$ in the definition of $\psi_3$ replaced by $4\times 10^{-7}$.
If $m_1\geq 10^{14}$, then we have $2.97+\log\psi_3(m_1)<-0.377493\cdots $ and \eqref{eq59} still holds
but this is incompatible with \eqref{eq57} when $m_2\geq \exp(294622)$.

If $10^7\leq m_1<10^{14}$ and $m_2\geq e^{40000}$, then we have $2.97+\log\psi_3(m_1)<-0.37723373\cdots $
and \eqref{eq59} still holds with $-0.3772337\cdots $ in place of $-0.377494$.
We proceed as just above and obtain $m_2<\exp(219207)$.
Similarly, if $5000\leq m_1<10^7$ and $m_2\geq e^{40000}$, then we have \eqref{eq59} with
$-0.363734\cdots $ in place of $-0.377494$ and $m_2<\exp(161128)$.
Similarly, if $53\leq m_1<5000$ and $m_2\geq e^{40000}$, then we have \eqref{eq59} with
$-0.234726\cdots $ in place of $-0.377494$ and $m_2<\exp(138880)$.

Moreover, if $5\leq m_1\leq 41$ and $m_2$ is at least the constant in the corresponding row in Table \ref{tbl6},
then we obtain \eqref{eq59} with $A_4$ taken from Table \ref{tbl4} with $\tilde m_2=e^{40000}$,
which is incompatible with \eqref{eq57}.

Now we are in turn to settle the remaining case $k_1=4$.
Then, we must have $x_1=239$ and $m_1=13$.
By Machin's formula, we have
\begin{equation}
4y_1 \arctan \frac{1}{5}\pm 4 \ell_3 \arctan \frac{1}{x_2} \pm 4 \ell_2 \arctan \frac{1}{x_3}=\frac{r^\prime\pi}{4}
\end{equation}
for some integer $r^\prime$, which may be zero.
If $r^\prime\neq 0$, then we have a reduced case $k_1=1$ with $x_1=13$ in place of $x_1=239$.
If $r^\prime=0$, then
\begin{equation}
y_1 \arctan \frac{1}{5}=\pm \ell_3 \arctan \frac{1}{x_2} \pm \ell_2 \arctan \frac{1}{x_3}.
\end{equation}
Hence, 
\begin{equation}
\abs{y_1} \arctan \frac{1}{5}<\frac{\ell_2+\ell_3}{x_2-1}\leq \frac{2G}{x_2-1}
\end{equation}
and therefore $m_2<((2G/\arctan(1/5))+1)^2+1$.
If $m_2>10^{10}/13$, then 
\begin{equation}
G\leq \floor{A_4(\log 13)\log^2(A_5\log 13)}=7701562.
\end{equation}
Hence, we have $m_2<6.089\times 10^{15}$.

Hence, we must have $m_1<2.531\times 10^{24}$ and $m_2<\exp(294622)$.
Moreover, if $m_1\geq 9.134\times 10^{22}$, then $m_2<e^{40000}$.
If $m_1\geq 10^{14}$ and $m_2\geq e^{40000}$, then, Theorem \ref{thm41} with
$(\tilde m_1, \tilde m_2)=(10^{14}, e^{40000})$ yields that
$\log x\leq (\log 2+k_i\log m_1+\ell_i\log m_2)/2<5.07\times 10^{12}$,
$x_i<\exp(5.07\times 10^{12})$,
and $\abs{y_i}\leq 2FG<1.316\times 10^{19}$.
Similarly, applying Theorem \ref{thm41} with $(\tilde m_1, \tilde m_2)$ according to
Table \ref{tbl6} in each case,
we can confirm that $x_i<\exp(5.07\times 10^{12})$ and $\abs{y_i}<1.319\times 10^{19}$
whenever $m_1\geq 53$.

\begin{table}
\caption{Upper bounds for $\log x_i$ and $y_i ~ (i=1, 2, 3)$ (Case II)}
\begin{center}
\begin{small}
\begin{tabular}{| c | c | c | c | c | c |}
\hline
$m_1$ & $m_2$ & $(\tilde m_1, \tilde m_2)$ & $\log m_2<$ & $\log x_i<$ & $y_i<$ \\
\hline
\addstackgap[.5\dimexpr 2pt \relax]{$\geq 10^{14}$} &
$\geq e^{40000}$ & $(10^{14}, e^{40000})$ & $294622$ &
$5.054\times 10^{12}$ & $1.312\times 10^{19}$ \\
$[10^7, 10^{14}]$ & $\geq e^{40000}$ & $(10^7, e^{40000})$ & $219207$ &
$2.982\times 10^{12}$ & $1.007\times 10^{19}$ \\
$[5000, 10^7]$ & $\geq e^{40000}$ & $(5000, e^{40000})$ & $161128$ &
$1.501\times 10^{12}$ & $6.935\times 10^{18}$ \\
$[53, 5000]$ & $\geq e^{40000}$ & $(53, e^{40000})$ & $138880$ &
$1.153\times 10^{12}$ & $8.986\times 10^{18}$ \\
$\geq 5000$ & $[e^{2000}, e^{40000}]$ & $(5000, e^{2000})$ & $40000$ & 
$1.475\times 10^{12}$ & $7.735\times 10^{18}$ \\
$[53, 5000]$ & $[e^{2000}, e^{40000}]$ & $(53, e^{2000})$ & $40000$ &
$1.664\times 10^{10}$ & $1.300\times 10^{17}$ \\
$\geq 53$ & $[m_1^*, e^{2000}]$ & $(53, m_1^*)$ & $2000$ & $1.830\times 10^{11}$ & $2.889\times 10^{18}$ \\

$41$ & $\geq e^{40000}$ & $(41, e^{40000})$ & $118023$ & $3.655\times 10^{11}$ & $2.438\times 10^{18}$ \\
$37$ & $\geq e^{40000}$ & $(37, e^{40000})$ & $108986$ & $1.339\times 10^{11}$ & $3.644\times 10^{17}$ \\
$29$ & $\geq e^{40000}$ & $(29, e^{40000})$ & $113396$ & $2.215\times 10^{11}$ & $1.028\times 10^{18}$ \\
$17$ & $\geq e^{40000}$ & $(17, e^{40000})$ & $110031$ & $1.493\times 10^{11}$ & $5.715\times 10^{17}$ \\
$13$ & $\geq e^{40000}$ & $(13, e^{40000})$ & $116234$ & $2.953\times 10^{11}$ & $2.339\times 10^{18}$ \\
$5$ & $\geq e^{40000}$ & $(5, e^{40000})$ & $115315$ & $2.169\times 10^{11}$ & $2.027\times 10^{18}$ \\
$41$ & $[10^{24}, e^{40000}]$ & $(41, 10^{24})$ & $40000$ & $1.498\times 10^{11}$ & $1.208\times 10^{18}$ \\
$37$ & $[10^{24}, e^{40000}]$ & $(37, 10^{24})$ & $40000$ & $6.245\times 10^{10}$ & $2.160\times 10^{17}$ \\
$29$ & $[10^{24}, e^{40000}]$ & $(29, 10^{24})$ & $40000$ & $9.574\times 10^{10}$ & $5.444\times 10^{17}$ \\
$17$ & $[10^{24}, e^{40000}]$ & $(17, 10^{24})$ & $40000$ & $6.701\times 10^{10}$ & $3.170\times 10^{17}$ \\
$13$ & $[10^{24}, e^{40000}]$ & $(13, 10^{24})$ & $40000$ & $1.218\times 10^{11}$ & $1.156\times 10^{18}$ \\
$5$ & $[10^{24}, e^{40000}]$ & $(5, 10^{24})$ & $40000$ & $8.975\times 10^{10}$ & $1.001\times 10^{18}$ \\
$41$ & $[10^{16}, 10^{24}]$ & $(41, 10^{16})$ & $24\log 10$ & $223204573$ & $1.943\times 10^{15}$ \\
$37$ & $[10^{16}, 10^{24}]$ & $(37, 10^{16})$ & $24\log 10$ & $94449003$ & $3.577\times 10^{14}$ \\
$29$ & $[10^{16}, 10^{24}]$ & $(29, 10^{16})$ & $24\log 10$ & $143666899$ & $8.874\times 10^{14}$ \\
$17$ & $[10^{16}, 10^{24}]$ & $(17, 10^{16})$ & $24\log 10$ & $100899200$ & $5.202\times 10^{14}$ \\
$13$ & $[10^{16}, 10^{24}]$ & $(13, 10^{16})$ & $24\log 10$ & $181014051$ & $1.850\times 10^{15}$ \\
$5$ & $[10^{16}, 10^{24}]$ & $(5, 10^{16})$ & $24\log 10$ & $133264185$ & $1.598\times 10^{14}$ \\
$41$ & $[243902441, 10^{16}]$ & $(41, 243902441)$ & $16\log 10$ & $137892418$ & $1.112\times 10^{15}$ \\
$37$ & $[270270277, 10^{16}]$ & $(37, 270270277)$ & $16\log 10$ & $57513230$ & $1.990\times 10^{14}$ \\
$29$ & $[344827589, 10^{16}]$ & $(29, 344827589)$ & $16\log 10$ & $88172538$ & $5.014\times 10^{14}$ \\
$17$ & $[588235297, 10^{16}]$ & $(17, 588235297)$ & $16\log 10$ & $61714770$ & $2.920\times 10^{14}$ \\
$13$ & $[769230773, 10^{16}]$ & $(13, 769230773)$ & $16\log 10$ & $112106482$ & $1.064\times 10^{15}$ \\
$5$ & $[2000000005, 10^{16}]$ & $(5, 2000000005)$ & $16\log 10$ & $82656260$ & $9.218\times 10^{14}$ \\
\hline
\end{tabular}
\label{tbl6}
\end{small}
\end{center}
\end{table}

If $m_1\leq 41$, then,
applying \eqref{eq51} with $\tilde m_2=p(10^{10}/m_1)$, $10^{16}$, $10^{24}$, and then $e^{40000}$, we have
$x_i<\exp(3.674\times 10^{11})$ and $\abs{y_i}\leq 2FG<2.451\times 10^{18}$ as in Table \ref{tbl6},
where $p(10^{10}/m_1)$ denotes the smallest integer above $10^{10}/m_1$ composed of primes $\equiv 1\Mod{4}$.

Thus we conclude that $x_i<\exp(5.07\times 10^{12})$, and $\abs{y_i}<1.316\times 10^{19}$ in Case II.

Finally, we assume that $m_1 m_2<10^{10}$.
It immediately follows from the latter part of Theorem \ref{thm41} that
\begin{equation}
k_i\log m_1+\ell_i\log m_2<1.7628 m_1 m_2\leq 1.7628(10^{10}-3),
\end{equation}
where we observe that $m_1 m_2\leq 10^{10}-3$ since $m_1$ and $m_2$ must be composed of primes $\equiv 1\Mod{4}$.
Hence, we obtain $\log x_i<(\log 2+k_i\log m_1+\ell_i\log m_2)/2<8813999998$
and, observing that $FG(\log m_1)(\log m_2)<(F\log m_1+G\log m_2)^2/4<(0.8814m_1 m_2)^2$,
\begin{equation}
\abs{y_i}\leq 2FG<\frac{2(0.8814m_1 m_2)^2}{\Omega}<\frac{2(8.814\times 10^9)^2}{\log 5\log(10^{10}/5)}<4.508\times 10^{18}.
\end{equation}

This completes the proof of the Theorem.

\section*{Acknowledgement}

This work was supported by the Research Institute for Mathematical Sciences,
a joint usage/research center in Kyoto University.

{}

\begin{thebibliography}{}
\bibitem{Cal}
Jack~S.~Calcut, Gaussian integers and arctangent identities for $\pi$,
\textit{Amer.~Math.~Monthly} \textbf{116} (2009), 515--530.

\bibitem{Car}
R.~D.~Carmichael, On the numerical factors of the arithmetic forms $\alpha^n\pm\beta^n$,
\textit{Ann.~Math.} \textbf{15} (1913), 30--70.

\bibitem{Eve}
J.-H.~Evertse, On equations in $S$-units and the Thue-Mahler equation,
\textit{Inv.~Math.} \textbf{75} (1984), 561--584.

\bibitem{HC1}
Henri~Cohen, \textit{Number Theory, Volume I: Tools and Diophantine Equations}, Springer, 2007.

\bibitem{HC2}
Henri~Cohen, \textit{Number Theory, Volume II: Analytic and Moderl Tools}, Springer, 2007.

\bibitem{Lau}
Michel~Laurent,
Linear forms in two logarithms and interpolation determinants II,
\textit{Acta~Arith.} \textbf{133} (2008), 325--348.

\bibitem{LMN}
Michel~Laurent, Maurice~Mignotte, and Yuri~Nesterenko,
Formes lin\'{e}aires en deux logarithmes et d\'{e}terminants d'interpolation,
\textit{J. Number Theory} \textbf{55} (1995), 285--321.

\bibitem{Leb}
M.~Lebesgue, Sur L'impossibilit\'{e}, en nombres entiers, de l'\'{e}quation $x^m=y^2+1$,
\textit{Nouv.~Ann.~Math.} \textbf{9} (1850), 178--180.

\bibitem{Lju}
W.~Ljunggren, Zur theorie der Gleichung $X^2+1=DY^4$,
\textit{Avh.~Norske,~Vid.~Akad.~Oslo} \textbf{1}, No.~5 (1942).

\bibitem{MV1}
Maurice~Mignotte and Paul~Voutier,
A kit on linear forms in three logarithms, with an appendix by Michel~Laurent,
\url{https://arxiv.org/abs/2205.08899v1}.

\bibitem{MV2}
Maurice~Mignotte and Paul~Voutier,
A kit on linear forms in three logarithms, with an appendix by Michel~Laurent,
\textit{Math. Comp.} \textbf{93} (2014), 1903--1951.

\bibitem{ST}
Ray~Steiner and Nikos~Tzanakis,
Simplifying the solution of Ljunggren's equation $X^2+1=2Y^4$,
{\it J.~Number~Theory} \textbf{37} (1991), 123--132.

\bibitem{Stm1}
Carl~St{\o}rmer, Solution compl\`{e}te en nombres entiers $m, n, x, y, k$ de l'\'{e}quation $m\textrm{arc tg}\frac{1}{x}+n\textrm{arc tg}\frac{1}{y}=k\frac{\pi}{4}$,
\textit{Skrift.~Vidensk.~Christiania I.~Math.-naturv.~Klasse} \textbf{1895}, Nr.~11, 21 pages.

\bibitem{Stm2}
Carl~St{\o}rmer, Sur l'application de la th\'{e}orie des nombres entiers complexes a la solution en nombres rationnels
$x_1, x_2, \ldots, x_n, c_1, c_2, \ldots, c_n, k$ de l'\'{e}quation: $c_1\textrm{arc tg} x_1+c_2\textrm{arc tg} x_2+\cdots +c_n\textrm{arc tg} x_n=k\frac{\pi}{4}$,
\textit{Arch.~Math.~Naturv.} \textbf{19} (1896), Nr.~3, 96 pages.

\bibitem{Stm3}
Carl~St{\o}rmer, Quelques th\'{e}or\`{e}mes sur l'\'{e}quation de Pell $x^2-Dy^2=\pm 1$ et leurs applications,
\textit{Skrift.~Vidensk.~Christiania~I.~Math.-naturv.~Klasse} \textbf{1897}, Nr.~2, 48 pages.

\bibitem{Twe}
Ian Tweddle, John Machin and Robert Simson on inverse-tangent series for $\pi$,
\textit{Arch.~Hist.~Exact~Sci.} \textbf{42} (1991), 1--14.

\bibitem{WH}
Michael~Roby~Wetherfield~and~Hwang~Chien-lih, 
Computing $\pi$, lists of Machin-type (inverse cotangent) identities for $\pi/4$,
\url{https://web.archive.org/web/20240204042153/https://machination.eclipse.co.uk/}.

\bibitem{Wol}
J.~Wolfskill, Bounding squares in second order recurrence sequences,
\textit{Acta~Arith.} \textbf{54} (1989), 127--145.

\bibitem{Wre}
J.~W.~Wrench~Jr., On the derivation of arctangent equalities,
Amer.~Math.~Monthly \textbf{45} (1938), 108--109.

\bibitem{Ybt}
Minoru~Yabuta, A simple proof of Carmichael's theorem on primitive divisors,
\textit{Fibonacci~Quart.} \textbf{39} (2001), 439--443.

\bibitem{Ymd}
Tomohiro~Yamada, A generalization of the Ramanujan-Nagell equation,
\textit{Glasgow~Math.~J.} \textbf{61} (2019), 535--544.

\bibitem{Ymd2}
Tomohiro~Yamada, A note on Laurent's paper on linear forms in two logarithms:
The argument of an algebraic power, \url{https://arxiv.org/abs/1906.00419}.

\end{thebibliography}
\end{document}